# ANALYTIC URNS

By Philippe Flajolet, Joaquim Gabarró and Helmut Pekari

*INRIA Rocquencourt, Universitat Politècnica de Catalunya and Universitat Politècnica de Catalunya*

This article describes a purely analytic approach to urn models of the generalized or extended Pólya–Eggenberger type, in the case of *two* types of balls and constant "balance," that is, constant row sum. The treatment starts from a quasilinear first-order partial differential equation associated with a combinatorial renormalization of the model and bases itself on elementary conformal mapping arguments coupled with singularity analysis techniques. Probabilistic consequences in the case of "subtractive" urns are new representations for the probability distribution of the urn's composition at any time $n$, structural information on the shape of moments of all orders, estimates of the speed of convergence to the Gaussian limit and an explicit determination of the associated large deviation function. In the general case, analytic solutions involve Abelian integrals over the Fermat curve $x^h + y^h = 1$. Several urn models, including a classical one associated with balanced trees (2–3 trees and fringe-balanced search trees) and related to a previous study of Panholzer and Prodinger, as well as all urns of balance 1 or 2 and a sporadic urn of balance 3, are shown to admit of explicit representations in terms of Weierstraß elliptic functions: these elliptic models appear precisely to correspond to regular tessellations of the Euclidean plane.

**0. Introduction.** In this study, we revisit the most basic urn model, namely the "generalized" (or "extended") Pólya–Eggenberger urn model with *two types of balls*, as described in the reference book of Johnson and Kotz (1977). Under this model an urn may contain two types of balls, say "black" (B) and "white" (W). The composition of the urn at time 0 is fixed. At time $n$, a ball in the urn is randomly chosen and its color is *observed* (thus the ball is selected, examined and then placed back into the urn): if it is black, then $\alpha$ black and $\beta$ white balls are subsequently inserted; if it is white, then $\gamma$ black balls and $\delta$ white balls are inserted. The evolution rule is then summarized by a $2 \times 2$-matrix

| drawn | added | |
|---|---|---|
| ↓ | B | W |
| B | $\alpha$ | $\beta$ |
| W | $\gamma$ | $\delta$ |

Negative values of the diagonal entries $\alpha, \delta$ are permissible and interpreted as an extraction (rather than an insertion) of balls; a model with both diagonal entries negative will be called here an *urn with subtraction* (of balls of the color chosen). The off-diagonal entries $\beta, \gamma$ are always taken to be nonnegative.

The urn model is said to be *balanced* if $\alpha + \beta = \gamma + \delta$, in which case the common sum of the matrix rows is the *balance*, denoted throughout by $s$. The $2 \times 2$ urn model may lead to widely differing behaviors depending on the values of the integer entries $\alpha, \beta, \gamma$ and $\delta$. For instance, Kotz, Mahmoud and Robert (2000) mention the (balanced) urn with matrix $\begin{pmatrix} 4 & 0 \\ 3 & 1 \end{pmatrix}$ for which the number of white balls picked in $n$ steps grows stochastically like $n^{1/4}$. Strikingly, Kotz, Mahmoud and Robert (2000) also study the (imbalanced) urn associated to $\begin{pmatrix} 1 & 0 \\ 1 & 1 \end{pmatrix}$ and show the corresponding number to be $\sim n/\log n$ in probability under a Poisson model. We do not address in this article models with more than two colors; see the paper of Smythe (1996) for a thorough probabilistic treatment, the works of Aldous (1991), Aldous, Flannery and Palacios (1988) for a discussion of almost sure convergence issues, and the comprehensive and independent recent studies of Janson (2004, 2005).

Our interest throughout this article is in urn models that are balanced. The conditions of having a matrix

$$(1) \qquad M = \begin{pmatrix} \alpha & \beta \\ \gamma & \delta \end{pmatrix} \qquad \text{with } \alpha + \beta = \gamma + \delta = s, \beta \geq 0, \gamma \geq 0,$$

are invariably assumed. We also allow ourselves on occasion to describe $M$ linearly as $(\alpha, \beta; \gamma, \delta)$. In such a case, each elementary action on the urn results in having the total number of balls increase by the fixed quantity $s$, so that the population at time $n$ has a predictable cardinality, which is exactly $t_0 + sn$ if $t_0$ is the initial size at time 0. For urns involving subtraction, certain simple arithmetic conditions on the parameters, called *tenability* (the Webster dictionary defines "tenable" as meaning "capable of being maintained"), ensure that the process cannot be "blocked"; these conditions are recalled in Section 1, (6), and are assumed to hold.

Balanced $2 \times 2$ urn models have been in particular considered by Bagchi and Pal (1985) who show the following: under a supplementary technical condition, namely that the ratio between eigenvalues of the matrix, $(\alpha-\beta)/(\alpha+\beta)$, lies in $(-\infty,\frac{1}{2})$, the distribution of the number of balls of one color obeys in the limit a normal distribution. Gouet (1993) further shows, under the assumptions of Bagchi and Pal, convergence of the discrete urn evolution to a stochastic Gaussian process, and he also investigates other cases using martingale arguments. Aldous, Flannery and Palacios (1988) observe that such results can be supplemented by almost sure convergence properties—their treatment extends the relation between branching processes and urn models to be found in the book by Athreya and Ney [(1972), Section V.9]. Thanks to the works of these and many other authors, the normal evolution of the process in the central regime can thus be regarded as well understood.

In this article, we revisit urn models under the radical angle of analysis. [Aldous (1991) otherwise provides an insightful comparison of the scopes of the traditional probabilistic approach and the modern methods of analysis of algorithms in his introductory section.] Our main results provide a complete analytic solution describing the composition of the urn at each instant, but, although our methods potentially apply to all the $2\times 2$ balanced urn schemes, we focus attention in this paper on urns involving *subtraction*, that is, having negative diagonal entries. The matrix can accordingly be taken under the form $(-a,a+s;b+s,-b)$, with balance $s\geq 0$ and diagonal coefficients $-a,-b<0$. (The urn's initial composition is fixed with $t_0$ balls in total of which $a_0$ are of the first type.) Such models with negative diagonal entries are occasionally mentioned by some authors as a harder nut to crack, since the direct embedding of urn schemes into branching processes explained in the book of Athreya and Ney [(1972), Section V.9] ceases to be directly applicable. [This position is perhaps to be taken with caution given the discussion in Aldous, Flannery and Palacios (1988) of extensions of the classical probabilistic framework.]

In the *first part* of the article (Sections 1 and 2), we introduce the partial differential equation approach to urn models with two types of balls and constant row sum. Our analysis starts with a *partial differential equation* (PDE) that is linear of the first order and that describes exactly snapshots of the urn compositions at all times. The solution of this partial differential equation, obtained by the standard method of "characteristics," provides an indirect expression for a bivariate generating function that encodes the possible configurations of the urn at each time $n$. It is found that this bivariate generating function is expressed in terms of a fundamental function $\psi$, which is defined implicitly by an equation of the form

$$\psi(I(u))=Q(u,v). \tag{2}$$

There $(u,v)$ lies on a Fermat curve $u^h+v^h=1$ with $h=a+b+s$ a sort of "complexity index," the quantity $Q$ being a rational function on the curve, and $I(u)$ an Abelian integral on that same curve—that is, the integral of a rational function on the curve. The parameterization (2) suffices in all cases to determine the dominant singularities of $\psi$ together with the associated singular expansions. As a consequence, analytic principles provide the (known) Gaussian law for the urn's composition at large instants, together with a precise determination of the speed of convergence as well as an explicit form of the large deviation function in terms of the Abelian integral $I(u)$ [see (3) for a specific instance]. In general, $\psi$ is associated with algebraic curves of genus strictly higher than 1. (Note: The Fermat curve is of high topological genus [Lang (1982)], namely $g = (h-1)(h-2)/2$, so that one already has $g = 10$ in the case of the $\mathcal{T}_{2,3}$ model discussed below. This makes the occasional existence of elliptic function solutions, which are objects of genus 1, quite remarkable.)

Our investigations were initially motivated by a desire to understand the specific urn model $\mathcal{T}_{2,3} := \begin{pmatrix} -2 & 3 \\ 4 & -3 \end{pmatrix}$, which forms the subject of the *second part* of this article. This particular urn process intervenes as a model of several schemes for managing an important data structure of computer science known as the search tree [Knuth (1998) and Mahmoud (1992)] and it surfaces in the analysis of 2–3 trees and fringe-balanced binary search trees [Aldous (1991), Aldous, Flannery and Palacios (1988), Bagchi and Pal (1985), Eisenbarth et al. (1982), Panholzer and Prodinger (1998) and Yao (1978)]. What is striking about this urn is that the model can be completely resolved in terms of *elliptic functions* of the Weierstraß type. For instance, our general results express that the probability of large deviations at time $n$ is exponentially small in $n$ with a rate that is a simple transform of the integral

$$K(u) := \frac{1}{(1-u^6)^{1/6}} \int_u^1 \frac{t}{(1-t^6)^{5/6}}\, dt. \tag{3}$$

A parallel elliptic connection had been uncovered earlier by Panholzer and Prodinger (1998) using rather different methods. Their penetrating analysis depends on the specific relationship that the $\mathcal{T}_{2,3}$ model entertains with a special type of "fringe-balanced" search trees—a root decomposition of the tree then leads to a perturbed nonlinear ordinary differential equation (of the rough form $Y''' = Y'^2 + \cdots$) akin to the one satisfied by the Weierstraß $\wp$-function. In this particular case, our elliptic connection for the $\mathcal{T}_{2,3}$ urn model could alternatively be deduced by reverse-engineering of the Panholzer–Prodinger treatment, combined with an easy reduction of a special urn model studied by Mahmoud (1998). We do not proceed along those lines since Panholzer and Prodinger's nonlinear differential approach is problem-specific, and, for example, it would not yield the other elliptic cases listed in (4).

The general character of our analytic results, Theorems 1 and 2, actually permits us to single out all the cases where elliptic solutions prevail, namely, all urns of balance 1 or 2 and a sporadic urn of balance 3, corresponding to the *six* matrices

$$(4)\qquad \begin{pmatrix}-2 & 3\\ 4 & -3\end{pmatrix},\quad \begin{pmatrix}-1 & 2\\ 3 & -2\end{pmatrix},\quad \begin{pmatrix}-1 & 2\\ 2 & -1\end{pmatrix},$$
$$\begin{pmatrix}-1 & 3\\ 3 & -1\end{pmatrix},\quad \begin{pmatrix}-1 & 3\\ 5 & -3\end{pmatrix},\quad \begin{pmatrix}-1 & 4\\ 5 & -2\end{pmatrix}.$$

We have chosen to illustrate the specificity of the elliptic models in a concrete way by developing properties of the $\mathcal{T}_{2,3}$ model (Section 3) and then proving our classification theorem regarding the elliptic cases (Theorem 4 and Section 4).

**1. Analytic solution of the general case.** We now take up the *general case* of a balanced urn model with two types of balls and negative diagonal entries. The matrix is of the form

$$(5)\qquad M = \begin{pmatrix}-a & a+s\\ b+s & -b\end{pmatrix}, \qquad a, b > 0,$$

with $s > 0$ the balance. Start with $a_0$ balls of the first type ("black") and $b_0$ balls of the second type ("white"), so that $t_0 = a_0 + b_0$ is the initial size; the size of the urn at time $n$ is then exactly the *deterministic* quantity $t_n = t_0 + ns$. In order for the urn not to be blocked by an infeasible request, the usual "tenability" conditions [Bagchi and Pal (1985) and Gouet (1993)] for urns with subtraction are assumed:

$$(6)\qquad \begin{cases}(\mathbf{T}_0): a \text{ divides } a_0 \text{ and } b \text{ divides } b_0;\\ (\mathbf{T}_1): a \text{ divides } b+s \text{ and } b \text{ divides } a+s.\end{cases}$$

We shall see soon that all such models are "solvable by quadrature" in the sense of Taylor [(1996), page 86]. In other words, only elementary algebraic functions, composition and inversion, as well as integration are involved in the solution, as is expressed by the general statement of Theorem 1. There results a complete characterization of dominant singularities, as summarized by Theorem 2. Probabilistic consequences are subsequently explored in Section 2.

1.1. *Algebraic approach.* Based on formal operator calculus, there is an elegant symbolic approach to the derivation of PDE's for urn models, which establishes a transparent connection between the combinatorial structure of a model and the PDE that expresses it.

The *combinatorial model* considers all balls involved in the game to be distinguished by distinct integer stamps: balls present at time 0 are stamped, say, $1, \dots, a_0$ for type B and $a_0 + 1, \dots, t_0$ for type W. New balls are stamped with "new" numbers: the balls that are taken away from the urn are (conventionally) the ball selected as well as others taken according to a deterministic policy, for example, by starting from smallest numbers. For instance, the urn $\begin{pmatrix} -1 & 2 \\ 3 & -2 \end{pmatrix}$ initialized with two balls of type B stamped with 1 and 2 may give rise to an evolution history starting as

| **Time** | **0** | **1** | **2** | **3** | |
|---|---|---|---|---|---|
| | choose 2 | choose 3 | choose 6 | choose 1 | |
| Urn | $\overbrace{1_{\rm B}, 2_{\rm B}}$, | $\overbrace{1_{\rm B}, 3_{\rm W}, 4_{\rm W}}$, | $\overbrace{1_{\rm B}, 5_{\rm B}, 6_{\rm B}, 7_{\rm B}}$, | $\overbrace{1_{\rm B}, 5_{\rm B}, 7_{\rm B}, 8_{\rm W}, 9_{\rm W}}$, | $\cdots$ |

with subscripts indicating colors/types of the corresponding balls.

In what follows, we consistently use $[z^n]f$ to denote the coefficient of $z^n$ in the formal power series or analytic function $f$.

One first needs to relate combinatorics and probability. We let $X_n$ be the random number of balls of type B at time $n$, and denote by $p_n(u)$ its *probability* generating function (PGF). Let $h_n(u)$ be the *counting* generating function of the evolution histories of length $n$, where $u$ marks the number of balls of type B: the coefficient $[u^k]h_n(u)$ is the number of histories comprising $n$ transformations of the urn and resulting in $k$ balls of type B. We have $p_0(u) = u^{a_0}$ as well as $h_0(u) = u^{a_0}$, and in general

$$p_n(u) = \frac{h_n(u)}{t_0(t_0+s)\cdots(t_0+(n-1)s)}, \tag{7}$$

since the total number of possible histories of length $n$ is

$$t_0(t_0+s)\cdots(t_0+(n-1)s) = n!s^n \binom{n+t_0/s-1}{n}, \tag{8}$$

as results from multiplication of $n$ elementary choices. (Naturally, the balance condition is crucial to this connection.) Introduce finally the exponential generating function of the $h_n(u)$, so that

$$H(z,u) := \sum_{n\geq 0} h_n(u)\frac{z^n}{n!} \tag{9}$$

is a *bivariate* generating function (BGF). As $u \to 1$, the bivariate generating function $H(z,u)$ degenerates into a simple algebraic function $H(z,1) = (1-sz)^{-t_0/s}$, since it then only counts histories in accordance with (8). Thus, $H(z,u)$ is a priori a "deformation" of a simple algebraic function.

For $u$ a variable, we let $\partial_u \equiv \frac{\partial}{\partial u}$ be the corresponding partial differential operator. It is notationally convenient to make use of the modified operator

$$\theta_u = u\partial_u \quad \text{so that} \quad \theta_u f = u\frac{\partial f}{\partial u}.$$

Differential operators are well known to correspond combinatorially to a "pointing" operation. For instance, one has

$$\partial_u u^a = au^{a-1}, \qquad \theta_u u^a = au^a,$$

so that $\partial_u$ may be interpreted as "select a $u$-element in all possible ways and remove it" while $\theta_u$ means "select a $u$-element in all possible ways and keep it." There are many instances in the combinatorics literature of such a usage of differential operators; see, for example, Bergeron, Labelle and Leroux [(1998), Section 2.1], Flajolet and Sedgewick (2003) and Goulden and Jackson [(1983), page 45].

Consider now an urn model defined by a matrix $M$ of the form (5), and represent momentarily a particular urn configuration with $\lambda$ white balls and $\mu$ black balls by the monomial $\mathfrak{m}_{\lambda,\mu} = u^\lambda v^\mu$. The partial differential operator (associated to $M$),

$$\Upsilon = u^{-a}v^{s+a}\theta_u + u^{s+b}v^{-b}\theta_v, \tag{10}$$

is such that the application of $\Upsilon$ to $\mathfrak{m}_{\lambda,\mu}$ describes all the possible successors of the urn represented by $\mathfrak{m}_{\lambda,\mu} = u^\lambda v^\mu$ when one step of ball replacement is performed.

Start with an urn of initial type $(a_0, b_0)$ represented by $u^{a_0}v^{b_0}$. Let $\widehat{h}_n(u, v)$ be the polynomial describing all possible evolutions of the urn in $n$ steps. [In particular, $\widehat{h}_n(u, 1) = h_n(u)$.] Then, one has

$$\widehat{h}_n(u, v) = \Upsilon^n \circ (u^{a_0}v^{b_0}).$$

We opt for exponential generating functions and define

$$\widehat{H}(z, u, v) = \sum_{n\geq 0} \widehat{h}_n(u, v)\frac{z^n}{n!}.$$

One has symbolically

$$\widehat{H}(z, u, v) = e^{z\Upsilon} \circ (u^{a_0}v^{b_0}),$$

where the exponential of operators is defined in the usual way:

$$e^{z\Upsilon} \circ g := \sum_{n\geq 0} \frac{z^n}{n!}(\Upsilon^n \circ g).$$

Then, the definition of the exponential immediately implies the differential relation

$$\partial_z(e^{z\Upsilon} \circ g) = \Upsilon e^{z\Upsilon} \circ g.$$

In other words, $\widehat{H}$ satisfies the PDE

$$\partial_z \widehat{H} = \Upsilon \circ \widehat{H}. \tag{11}$$

The last equation is almost the PDE we are looking for but not quite (it has a supplementary variable, $v$). Given the balance condition, the urn population increases by exactly $s$ at each step. Accordingly, $\widehat{H}$ involves three variables, $u, v$ and $z$, but their exponents in $\widehat{H}$ are bound by a homogeneity condition, each monomial generated being of the form $u^\lambda v^\mu z^n$ with $\lambda+\mu = sn+t_0$. In other words, each monomial $\mathfrak{m}$ composing $H$ satisfies

$$(\theta_u + \theta_v - s\theta_z)\mathfrak{m} = t_0\mathfrak{m}, \tag{12}$$

and the relation extends by linearity to $\widehat{H}$ itself.

In summary, a system of two equations now determines $\widehat{H}$ (with $\theta_u \equiv u\partial_u$):

$$\begin{aligned} \partial_z \widehat{H} &= \Upsilon \circ \widehat{H}, \\ (\theta_u + \theta_v - s\theta_z)\widehat{H} &= t_0 \widehat{H}. \end{aligned} \tag{13}$$

One can then eliminate the explicit differential dependency on $v$ (the operator $\partial_v$) and get from (10) and (13)

$$\partial_z \widehat{H} = u^{-a} v^{1+a} \theta_u \widehat{H} + u^{1+b} v^{1-b} (s\theta_z \widehat{H} - \theta_u \widehat{H} - t_0 \widehat{H}).$$

At this stage it becomes possible to set $v = 1$, that is, completely eliminate the redundant variable $v$ itself. In this way one obtains the *fundamental PDE*

$$(1 - szu^{b+s})\frac{\partial H}{\partial z} + (u^{b+s+1} - u^{1-a})\frac{\partial H}{\partial u} - t_0 u^{b+s} H = 0, \tag{14}$$

where $H \equiv H(z,u) = \widehat{H}(z;u,1)$.

The main result is then:

THEOREM 1. *Consider the urn specified by matrix* $\left(\begin{smallmatrix} -a & a+s \\ b+s & -b \end{smallmatrix}\right)$, *with initial conditions* $(a_0, b_0)$ *and* $t_0 := a_0 + b_0$, *assuming it to be tenable. The probability generating function at time* $n$ *of the urn's composition is*

$$p_n(u) = \frac{\Gamma(n+1)\Gamma(t_0/s)}{s^n \Gamma(n + t_0/s)} [z^n] H(z,u),$$

*where the bivariate generating function* $H(z,u)$ *is given by*

$$H(z,u) = \delta(u)^{t_0} \psi(z\delta(u)^s + I(u)),$$

*with*

$$\delta(u) := (1-u^h)^{1/h}, \qquad I(u) := \int_0^u \frac{t^{a-1}}{\delta(t)^{a+b}}\,dt, \qquad h := a+b+s,$$

*and the function* $\psi$ *is defined implicitly by*

$$\psi(I(u)) = \frac{u^{a_0}}{\delta(u)^{t_0}}.$$

PROOF. We make use of the classical method of characteristics exposed in most textbooks, for example, Zwillinger [(1989), Section 94]. Following this method, one first associates to the linear first-order partial differential equation (14) the *ordinary differential system*

$$\frac{dz}{1 - szu^{b+s}} = \frac{du}{u^{s+b+1} - u^{1-a}} = \frac{dw}{t_0 u^{b+s} w}, \tag{15}$$

where $w$ "represents" $H$, and look for its first integrals.

The equation binding $w$ and $u$ allows for separation of variables,

$$\frac{dw}{w} = t_0 \frac{u^{h-1}}{u^h - 1}\, du,$$

so that a first integral of (15) is

$$w\delta(u)^{-t_0} = C_1. \tag{16}$$

The equation binding $z$ and $u$ is similar but inhomogeneous:

$$\frac{dz}{du} = -sz\frac{u^{h-1}}{u^h - 1} + \frac{u^{a-1}}{u^h - 1}.$$

The homogeneous equation is solved by separation of variables as $z = \xi \cdot (1 - u^h)^{-s/h}$. By the variation-of-constant technique, one finds

$$z = \xi(u)(1 - u^h)^{-s/h}, \qquad \xi(u) = -\int^u \frac{t^{a-1}}{(1 - t^h)^{(a+b)/h}}\, dt,$$

so that a first integral of (15) is

$$z\delta(u)^s + I(u) = C_2. \tag{17}$$

According to the method of characteristics, the general solution to the fundamental PDE (14) is obtained by coupling the two first integrals (16) and (17), namely

$$\Phi(H(z,u)\delta(u)^{-t_0}, z\delta(u)^s + I(u)) = 0,$$

for an arbitrary bivariate function $\Phi$. Solving symbolically for $H$ puts the solution in the form

$$H(z,u) = \delta(u)^{t_0}\psi(z\delta(u)^s + I(u)), \tag{18}$$

for an arbitrary univariate function $\psi$. The initial condition $H(0,u) = u^{a_0}$ finally identifies $\psi$ as defined implicitly through inversion of $I(u)$, namely, $\psi(I(u)) = u^{a_0}/\delta(u)^{t_0}$.

We observe next that $\psi(z)$ is analytic at 0. Indeed the tenability conditions of (6) imply that $a$ must divide $a_0$ and $a$ must divide $b + s$, hence $a$ divides $h = a + b + s$. In particular, the general form of the parameterization of $\psi$ near 0 is $\psi(u^a) \asymp u^{a_0}$, that is, $\psi(z) \asymp u^{a_0/a}$, which is compatible with analyticity. In fact, the expansions involved are of the form

$$\psi\left(u^a \sum_{j\geq 0} \lambda_j u^{jh}\right) = u^{a_0}\left(\sum_{j\geq 0} \mu_j u^{jh}\right),$$

for some real coefficients $\lambda_j, \mu_j$ and $u$ ranging in a small enough complex neighborhood of 0. Examination of the exponents involved in the inversion shows that $\psi(z)$ can be expanded as a power series in $z$, and analyticity of $\psi$ at 0 results. □

*Sensitivity to initial conditions.* When the initial state of the urn is changed, the functions involved still live in the same general class. Indeed, the $\psi$ function corresponding to an initial urn of composition $(a_0, b_0)$ factorizes, in accordance with Theorem 1, as

$$\psi(z) = \psi_I(z)^{a_0/a} \cdot \psi_{II}(z)^{b_0/b}, \tag{19}$$

where $\psi_I, \psi_{II}$ are determined implicitly by

$$\psi_I(I(u)) = \left(\frac{u}{\delta(u)}\right)^a, \qquad \psi_{II}(I(u)) = \left(\frac{1}{\delta(u)}\right)^b, \tag{20}$$

corresponding to an urn initialized with $t_0 = a_0 = a$ and $t_0 = b_0 = b$, respectively. The analytic treatment given below extends to both functions $\psi_I, \psi_{II}$, and it is seen that the main determinant of the category of special functions encountered is the index $h$ of the Fermat curve and the integral $I(u)$. Equations (19) and (20) thus give us flexibility for the choice of the initial conditions, as is done repeatedly below.

In the case where $a$ and $b$ are each at least $-1$, balls have a "descendance" and the evolution of descendants are combinatorially independent. Accordingly, the factorization (20) can be viewed as expressing the fact that the histories of all the initial balls can be freely shuffled. (It is known that shuffle products correspond to products of exponential generating functions.) A parallel decomposition underlies the probabilistic reduction of this class of urn models to multitype branching processes [Athreya and Ney (1972)], at least in the case where no diagonal entry is below $-1$, so that the disappearances of balls are not coupled.

1.2. *Complex-analytic structures.* For notational simplicity, we shall adopt in this section the initial conditions $a_0 = t_0 = a$, that is, the urn is initialized with exactly $a$ balls of the first type (B): by (19), (20) and the ensuing remarks, no essential loss of generality is implied by such a choice.

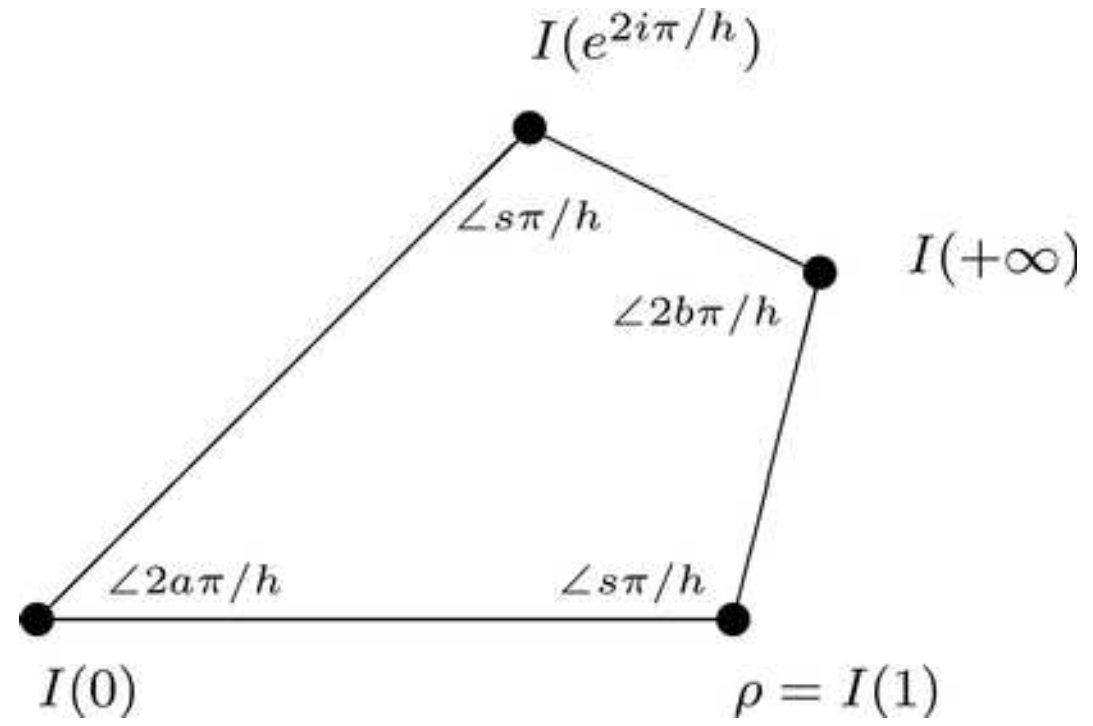

FIG. 1. *The elementary kite is the image of a small sector $S_0$.*

We make use of the quantity $h = a + b + s$. The function $v = \delta(u)$ corresponds to the complex Fermat curve,

$$u^h + v^h = 1,$$

which has topological genus $g = (h-1)(h-2)/2$. Following a classical terminology, the integral $I(u) \equiv \int u^{a-1} v^{-a-b}$ is a particular Abelian integral over this curve. The diagram that summarizes the parameterization of $\psi$ is then

$$\begin{array}{ccc} & u & \\ I(u) \swarrow & & \searrow J(u) \\ z & \xrightarrow{\psi} & \psi(z) \end{array} \qquad\qquad J(u) \equiv \frac{u^a}{\delta(u)^a}.$$

The major characteristics of an urn model turn out to be determined by the nature of the map $u \mapsto I(u)$ in the complex plane, with $J(u)$ playing only a secondary role.

As observed in the proof of Theorem 1, the function $\psi$ is analytic at 0 and it satisfies $\psi(z) \asymp z^{a_0/a}$ there. Also, the nature of the parameterization near 0, where $I(u) \asymp u^a$ implies that $I(u)$ effects an $a$-fold covering of a neighborhood of the origin and that $\psi(z)$ is of the form

$$\psi(z) = z^{a_0/a} \widehat{\psi}(z^{h/a}), \qquad (21)$$

for some $\widehat{\psi}$ analytic at the origin. In other words, in order to define $\psi$ parametrically by means of $u$, it suffices to let $u$ range in a sector $\mathcal{H}$ of angle $2\pi/(h/a)$ at the origin, and from now on, we shall do so. (As already noted, the tenability conditions precisely imply that $a$ divides $h$.)

Consider first the complex plane with $h$ rays emanating from 0 and having directions given by all the $h$th roots of unity. The sector $S_j$ is defined as

$$S_j := \left\{ z, z = Re^{i\theta}, 0 < R < \infty, \frac{2j\pi}{h} < \theta < \frac{2(j+1)\pi}{h} \right\}.$$

We claim (and prove below) that the image of $S_0$ by $I(u)$ is the interior of a quadrilateral (Figure 1), with vertices at the points

$$0, \qquad I(1), \qquad I(+\infty), \qquad I(e^{2i\pi/h}),$$

and call this quadrilateral $\mathcal{K}$ the *elementary kite*. [Note: The incomplete Beta integrals that make up $\psi$ are related to hypergeometric functions as well as to the Schwarz–Christoffel integrals of conformal mapping theory. For the latter aspects, see, e.g., the book of Nehari (1975), his Exercise 4, page 196, and his Chapter V.] One has

$$I(1) = \int_0^1 \frac{t^{a-1}}{(1-t^h)^{(a+b)/h}}\,dt = \frac{1}{h}B\left(\frac{a}{h}, \frac{s}{h}\right) = \frac{1}{h}\frac{\Gamma(a/h)\Gamma(s/h)}{\Gamma((a+s)/h)},$$

where use has been made of the usual Eulerian Beta integral [Whittaker and Watson (1927)]

$$B(\alpha,\beta) := \int_0^1 t^{\alpha-1}(1-t)^{\beta-1}\,dt = \frac{\Gamma(\alpha)\Gamma(\beta)}{\Gamma(\alpha+\beta)}. \tag{22}$$

We henceforth denote the quantity $I(1)$ by $\rho$.

The local mapping properties corresponding to the four vertices of the elementary kite are determined by the local behavior of $I(u)$: (i) at 0, $I(u)$ multiplies angles by $a$, so that the angle of the kite at 0 is $\frac{2\pi a}{h}$; (ii) at 1, $I(u)$ multiplies angles by $\frac{s}{h}$, so that the angle of the kite at vertex $I(1)$ is $\frac{\pi s}{h}$ [and similarly for vertex $I(e^{2i\pi/h})$]; (iii) at infinity, $I(u)$ multiplies angles by $b$, so that the angle at $I(+\infty)$ is $\frac{2\pi b}{h}$. In order to see that $I(u)$ maps the boundary of $S_0$ to that of $\mathcal{K}$, observe first that $I(u)$ maps $[0,1]$ onto the segment $[0, I(1)]$ by monotonicity of the integrand. Then, as $u$ continues to increase along $(1,+\infty)$ passing above 1, the function $\delta(u)$ becomes a complex number of fixed argument $\frac{a+b}{h}\pi$. In other words, $I$ maps the ray $[1,+\infty]$ to the segment $[I(1), I(+\infty)]$ (with $+\infty$ here understood as lying inside $S_0$). A similar discussion gives the mapping properties associated with the other two sides of the kite $\mathcal{K}$.

Next, we turn to sectors $S_1, \ldots$. Let $\zeta := e^{2i\pi/h}$. The image of sector $S_j$ is simply obtained as the image of $S_0$ by $I(u\zeta^j)$, which, by a linear change of variables [since $I(u\zeta^j) = \zeta^{-ja}I(u)$] is the image of the elementary kite under a rotation of angle $-\frac{2ja\pi}{h}$; see Figure 2 for a particular instance. Because of (21) and the accompanying remarks, it is sufficient to consider $0 \le j < \frac{h}{a}$.

Definition 1. The *fundamental polygon* of an urn model is the (closure of) the union of $h/a$ regularly rotated versions of the elementary kite about the origin.

We state:

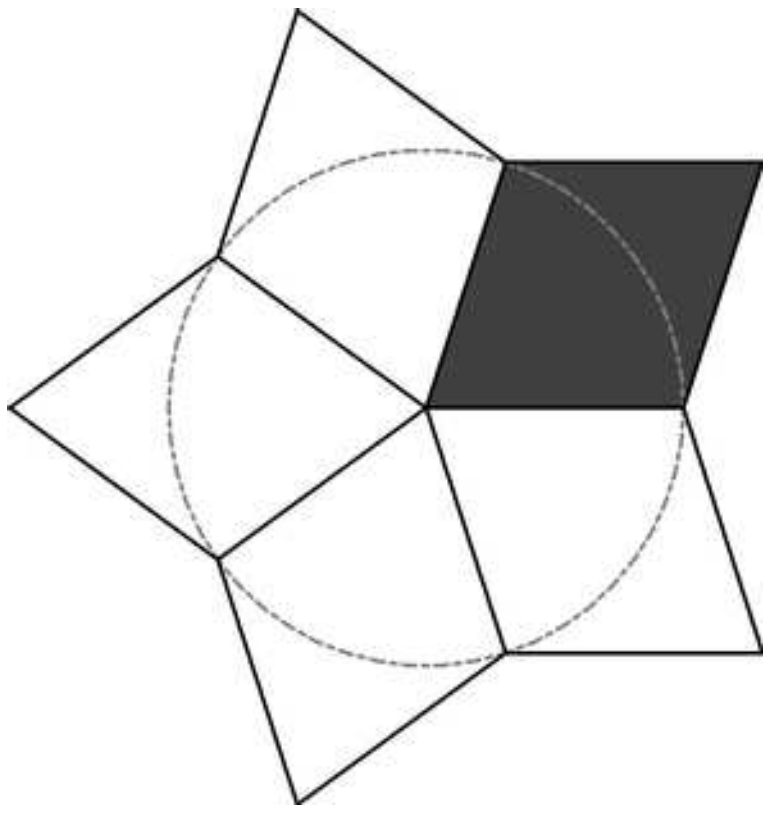

FIG. 2. *The elementary kite (in black) and the fundamental polygon associated with the urn* $(-1,4;4,-1)$.

THEOREM 2. *Consider a balanced* $2\times 2$ *urn with subtraction as in Theorem 1 and let it be initialized with* $a_0=a$, $b_0=0$. *The corresponding function* $\psi$ *is analytic for* $z$ *in the fundamental polygon of Definition 1. Furthermore, it is analytic in* $|z|<\rho$, *where*

$$\rho = I(1) = \int_0^1 \frac{t^{a-1}}{(1-t^h)^{(a+b)/h}}\,dt = \frac{1}{h}B\left(\frac{a}{h},\frac{s}{h}\right) = \frac{1}{h}\frac{\Gamma(a/h)\Gamma(s/h)}{\Gamma((a+s)/h)}.$$

*On* $|z|=\rho$, *the function* $\psi$ *is singular at* $\rho$ *and at the points* $\rho\omega^j$ *where* $\omega=\exp(2i\pi\frac{a}{h})$ *is an* $(h/a)$*th root of unity, regular at the other points. Its singular expansion as* $z\to\rho$ *is of the form*

$$\psi(z) = s^{-t_0/s}(\rho-z)^{-t_0/s}\mathcal{A}((\rho-z)^{h/s}), \tag{23}$$

*with* $\mathcal{A}$ *analytic at* 0, $\mathcal{A}(0)=1$, $\mathcal{A}'(0)\neq 0$. (*Principal determinations as* $z\to\rho^-$ *are assumed.*) *This expansion extends to a sector of opening larger than* $\pi$ *at* $\rho$.

*At the points* $z=\rho\omega^j$, *the singular expansion is determined from the expansion at* $z=\rho$ *by the fact that* $\psi(z)z^{-a_0/a}$ *is invariant under the mapping* $z\mapsto\omega z$.

The expansion (23) gives $\psi(z)$ as the product of a main singular part of the form $(\rho-z)^{-t_0/s}$ multiplied by a Puiseux series, that is, a series in fractional powers of $(\rho-z)$. We shall occasionally refer to the quantity $h/s$ as the *Puiseux exponent* of $\psi$. It plays a special role in the discussion of elliptic urns in Sections 3 and 4, in which case it reduces to an integer value.

PROOF OF THEOREM 2. First, the fact that $I(u)$ assumes each value in $\mathcal{K}$ once and only once when $u\in S_0$ is a consequence of basic properties of conformal mapping theory, which we recall. Let $\beta$ be an arbitrary number interior to $\mathcal{K}$. The number $\nu(\beta)$ of times that $I(u)$ assumes the value $\beta \in \mathcal{K}$ for $u$ interior to $S_0$ is by the residue theorem

$$\nu(\beta) = \frac{1}{2i\pi} \int_{\partial S_0} \frac{I'(u)}{I(u)-\beta}\, du,$$

where $\partial X$ represents the boundary of a region $X$ oriented positively. Then, the change of variables $I(u) = x$ gives

$$\nu(\beta) = \frac{1}{2i\pi} \int_{\partial \mathcal{K}} \frac{dx}{x-\beta} = 1,$$

where the reduction to the value 1 is due to the fact that $\beta$ is by assumption interior to $\mathcal{K}$. This implies that the functional inverse $u = I^{(-1)}(z)$ is well defined (and analytic) for $z$ interior to $\mathcal{K}$, and so is $\psi(z)$ since $\psi(z) = J(u)$ while $J(u)$ depends analytically on $u$. These properties extend in turn to the fundamental polygon by rotations of the base sector.

We next examine the behavior of $\psi$ near $\rho = I(1)$, corresponding to $u$ in the vicinity of 1 (say, $u \to 1^-$, to fix ideas). The expansion can be constructed by means of a local uniformizing parameter, here, $1-u = \tau^h$. Write

$$\delta(y) = \Delta(y)(1-y)^{1/h},$$

so that $\Delta(y)$ is analytic at $y = 1$. By the change of variables $u \mapsto 1-\tau^h$, one finds

$$\begin{aligned} I(1) - I(u) &= h \int_0^{\tau} (1-y^h)^{a-1} \Delta(1-y^h)^{-a-b} y^{s-1}\, dy, \\ &= \frac{1}{s}(h^{1/h}\tau)^s \left(1 + \frac{(h(b-a+2)-a-b)s}{2h(s+1)} \tau^h + \cdots \right), \\ J(u) &= \frac{(1-\tau^h)^{a_0}}{\Delta(1-\tau^h)^{t_0}} \\ &= (h^{1/h}\tau)^{-t_0} \left(1 - \frac{h(2a_0-t_0)+t_0}{2h} \tau^h + \cdots \right), \end{aligned} \tag{24}$$

where now $\tau \to 0$ corresponds to $u \to 1$ (the series expansions proceed by powers of $\tau^h$). Thus the parameterization is of the form

$$\rho - z = \frac{1}{s}(h^{1/h}\tau)^s U(\tau^h), \qquad \psi(z) = (h^{1/h}\tau)^{-t_0} V(\tau^h),$$

where $U, V$ are analytic at 0 and $U(0) = V(0) = 1$. By analytic inversion, this shows that there exists a full expansion of the type (23), with $\mathcal{A}$ analytic at 0. In other words, the point $\rho$ is a singularity of $\psi$ that is a branch point with dominant singular exponent equal to $-t_0/h$.

By rotational symmetry, an expansion of a nature similar to (23) also holds at the conjugate points $\rho\omega^j$ where $\omega = e^{2i\pi/h}$ is an $h$th root of unity.

Since $\psi(z)$ has nonnegative coefficients, it satisfies Pringsheim's theorem and is thus analytic in $|z| < \rho$. By the triangle inequality, we have $|I(ue^{i\theta})| \le I(u)$ for $u \in (0,1)$ and $\theta \in (-\pi, \pi)$. Since the nonzero terms composing the Taylor expansion of $I$ at the origin are of the form $u^{a+jh}$, the inequality $|I(ue^{i\theta})| < I(u)$ is strict as soon as $\theta$ is not a multiple of $\frac{2\pi}{h}$ and $I(u)$ is invertible. From there, it results that $\psi$ is analytic on $|z| = \rho$ except for the regularly spaced singularities quoted in the statement.

This provides the analytic continuation of $\psi$ in the fundamental polygon as well as in the disc of radius $\rho$. If the fundamental polygon is such that $s/h > \frac{1}{2}$, analytic continuation of $\psi$ outside its disc of convergence is granted and the proof is completed—this is the situation exemplified by Figure 2. Otherwise, the convergent character of $\mathcal{A}$ provides the analytic continuation of $\psi$ in sectors rooted at singularities and extending beyond the disc $z < \rho$. (This last situation is encountered in the $\mathcal{T}_{2,3}$ model detailed below; see Figures 3 and 4.) □

For instance, the urn $\binom{-1\ \ \ 4}{4\ \ -1}$ gives rise to the fundamental polygon displayed in Figure 2. One has $s = 3$, $h = 5$ and $\delta(u) = (1 - u^5)^{1/5}$, so that the fundamental polygon is a star with five branches. At the origin, we find $I(u) = u + \frac{1}{15}u^6 + \cdots$ and $\psi(z) = z + \frac{2}{15}z^6 + \cdots$. There is an algebraic branch point at $\rho$ where $\psi(\rho - x) \asymp (\rho - x)^{-1/3}$ and at the conjugate points $\rho\omega^j$ where $\omega^5 = 1$. The nature of the branch point of $\psi$ at $\rho$ is

$$\psi(z) = (3Z)^{-1/3}\left(1 - \tfrac{9}{40}(3Z)^{5/3} - \tfrac{1143}{10400}(3Z)^{10/3} + \cdots\right), \qquad Z := (\rho - z),$$

the Puiseux exponent associated with $\mathcal{A}$ being the fractional number $h/s = 5/3$.

**2. Probabilistic consequences.** Singularity analysis [Flajolet and Odlyzko (1990) and Odlyzko (1995)] makes it possible to extract very precise information on coefficients of a generating function once the function is recognized to have isolated singularities on the boundary of its disc of convergence. Under such conditions, assuming first unicity of the dominant singularity $\sigma$, an asymptotic estimate for a function $F$ of the form

$$F(z) \underset{z \to \sigma}{\sim} c(1 - z/\sigma)^{-\alpha}$$

valid in a complex region beyond $\sigma$ (a sector centered at $\sigma$ of opening angle larger than $\pi$ and including the disc of convergence) entails a matching estimate for the function's coefficients:

$$[z^n]F(z) \underset{n \to \infty}{\sim} c\sigma^{-n}\frac{n^{\alpha-1}}{\Gamma(\alpha)}.$$

Full asymptotic expansions can be transferred from functions to coefficients in a similar way [Flajolet and Odlyzko (1990)]. Also, in the presence of several dominant singularities, contributions to coefficients are to be composed additively. (This technology based on Hankel contours is of the complex Tauberian type.) It applies to the function $\psi(z)$ itself, and since $\psi(z) = H(z,0)$, it provides immediately sharp estimates of the probabilities that all balls are of the same color at epoch $n$, which corresponds to extreme large deviations.

COROLLARY 1 (Extreme large deviations). *For any balanced* $2 \times 2$ *urn with subtraction* (*i.e., negative diagonal entries*), *the probability that balls at time* $n$ *are all of the same color and of the second type* (W) *is*

$$\frac{h}{a}(s\rho)^{-n-t_0/s}\left(1+O\left(\frac{1}{n^{h/s}}\right)\right) \qquad \text{for } n \equiv \frac{a_0}{a}\left(\bmod \frac{h}{a}\right).$$

PROOF. The singularity at $z=\rho$ of $H(z,0)=\psi(z)$ contributes to $[z^n]\psi(z)$ a term

$$(s\rho)^{-t_0/s}\rho^{-n}\frac{n^{t_0/s-1}}{\Gamma(t_0/s)}\left(1+O\left(\frac{1}{n^{h/s}}\right)\right).$$

We have the periodicity expressed by (21). Thus, for $n$ in a suitable congruence class, there are $h/a$ similarly behaving singularities to be combined. The total number of histories of length $n$ is, from (8), asymptotic to

$$n!s^n\frac{n^{t_0/s-1}}{\Gamma(t_0/s)}.$$

The result follows after normalization by the latter quantity. $\square$

Next, we summarize the basic technology used to derive a Gaussian limit by the following statement, a simplified form of what is often referred to as the "quasi-powers theorem" and originates in works of Bender (1973) and Hwang (1998). Throughout this article, we use $\mathbb{E}$ and $\mathbb{V}$ to denote the expectation and variance operators.

LEMMA 1 (Quasi-powers theorem). *Let* $q_n(u) = \mathbb{E}(u^{Y_n})$ *be a family of probability-generating functions relative to discrete random variables* $Y_n$. *Assume that there exist two functions* $A(u), B(u)$ *analytic in a neighborhood* $\mathcal{V}$ *of* $u=1$, *such that, in this neighborhood the quasi-power approximation*

$$q_n(u) = A(u)B(u)^n(1+\varepsilon_n(u)) \qquad \text{as } n\to\infty \tag{25}$$

*holds, where* $|\varepsilon_n(u)| = O(n^{-1/2})$ *uniformly with respect to* $u$, *that is,* $\sup_{u\in\mathcal{V}}|\varepsilon_n(u)| = O(n^{-1/2})$. *Assume also the* variability condition

$$\sigma^2 \neq 0 \qquad \text{where } \sigma^2 := \lim_{n\to\infty}\frac{\mathbb{V}Y_n}{n} \equiv B''(1)+B'(1)-B'(1)^2. \tag{26}$$

[*Equivalence between the two forms of* $\sigma^2$ *is granted under condition* (25).] *Then, the random variables* $Y_n$ *converge in law to a Gaussian limit, with speed of convergence* $O(n^{-1/2})$*: for any* $x$*, one has*

$$\mathbb{P}\left(\frac{Y_n - \mathbb{E}(Y_n)}{\sqrt{\mathbb{V}Y_n}} \le x\right) = \frac{1}{\sqrt{2\pi}} \int_{-\infty}^{x} e^{-y^2/2}\, dy + O\left(\frac{1}{\sqrt{n}}\right). \tag{27}$$

Proof [Sketch; see Bender (1973) and Hwang (1998) for details]. The characteristic function $q_n(e^{it})$ of the $Y_n$ is by assumption closely approximated by an $n$th power. The variable $Y_n$ is next centered around its mean and scaled by its standard deviation in the usual way. A calculation similar to the usual case of independent random variables [e.g., Billingsley (1986), page 367] then shows the standardized version of $q_n(e^{it})$ to converge to $e^{-t^2/2}$, which is the characteristic function of a Gaussian law. The speed of convergence estimate finally results from the Berry–Esseen inequality found in Lukacs (1970). □

In order to apply the quasi-powers theorem, we choose a small complex neighborhood of $u=1$ and keep $u$ in this neighborhood. The BGF $H(z,u)$ rewrites as

$$H(z,u) = \delta(u)^{t_0} \psi(\rho - \delta(u)^s (K(u) - z)), \tag{28}$$

where

$$K(u) := \frac{1}{\delta(u)^s} \int_u^1 \frac{t^{a-1}}{\delta(u)^{a+b}}\, dt, \tag{29}$$

and $K(u)$ has a removable singularity at 1 with $K(1) = 1/s$. Treating $u$ as a parameter, we find that, as a function of $z$, the quantity $H(z,u)$ has a singularity at $z = K(u)$ that gets smoothly displaced when $u$ varies. Because of the nature of the singularity of $\psi$ at $\rho$, the singular exponent remains equal to the constant $-t_0/s$. Thus, for some function $L(u)$ that is analytic at $u=1$, one has

$$[z^n]H(z,u) = L(u) K(u)^{-n} n^{t_0/s-1} \left(1 + O\left(\frac{1}{n^{h/s}}\right)\right),$$

the error term being uniform by virtue of uniformity of the singularity analysis process [Flajolet and Odlyzko (1990)]. This has the shape of a bona fide quasi-powers approximation for the probability generating function,

$$p_n(u) = \frac{[z^n]H(z,u)}{[z^n]H(z,1)} = \frac{L(u)}{L(1)} \left(\frac{K(u)}{K(1)}\right)^{-n} \left(1 + O\left(\frac{1}{n^{h/s}}\right)\right).$$

The quasi-powers theorem then applies and gives:

COROLLARY 2 (Gaussian law and speed). *For any balanced* $2 \times 2$ *urn with subtraction, the random variable* $X_n$ *representing the number of balls of the first color* (B) *at time* $n$ *is asymptotically Gaussian with speed of convergence to the limit* $O(n^{-1/2})$, *as expressed by* (27).

The fact that the limit distribution is Gaussian was first observed by Bagchi and Pal (1985). These authors applied the moment method and determined the main asymptotic orders of moments of the centered variable $X_n - \mathbb{E}(X_n)$. Their method does not, however, appear to give access to the speed of convergence as expressed above. This speed is on the other hand neatly implied by the functional limit theorem of Gouet (1993). Here, we emphasize that the speed of convergence comes out almost immediately from the analytic approach.

In general, the moments are computable systematically from the exact expressions of Theorem 1 by successive differentiation with respect to $u$ upon setting $u = 1$ and making use of the singularities of $\psi$ and its derivatives as expressed by Theorem 2. All moments happen to be expressible in *closed form*.

COROLLARY 3 (Moments). *For any balanced* $2 \times 2$ *urn with subtraction and any* $r \geq 0$, *the* $r$*th factorial moment of the distribution of* $X_n$ *is of hypergeometric type: it is a finite linear combination of terms of the form*

$$\frac{\binom{n+t_0/s+\ell-kh/s-1}{n}}{\binom{n+t_0/s-1}{n}}, \qquad 0 \leq k, \ell \leq r.$$

The existence of such *finite* binomial forms for moments of *all* orders does not seem to have been previously noticed. Explicit forms are given by Kotz, Mahmoud and Robert (2000), but only for the first moment and at the cost of some labor, in the case of the urn model $(4, 0; 3, 1)$. Bagchi and Pal (1985) obtained such expressions for a wide class of urns, but in the case of the first two moments only.

PROOF OF COROLLARY 3. The quantity

$$\chi_r(z) := \frac{1}{r!}(\partial_u^r H(z,u))_{u=1} = [(u-1)^r] H(z,u)$$

is a generating function of the $r$th factorial moment of $X_n$ in the sense that

$$\mathbb{E}(X_n^{\underline{r}}) = \frac{[z^n]\chi_r(z)}{[z^n]\chi_0(z)},$$

where the notation $X^{\underline{r}}$ is the usual notation for falling factorials [Graham, Knuth and Patashnik (1989)], namely,

$$a^{\underline{r}} = a(a-1)\cdots(a-r+1). \tag{30}$$

In order to gain access to such moments, we make use of the singular expansion of $\psi$ near $\rho$. Given the variant form of $H(z,u)$ from (28), the singular expansion (23) of Theorem 2 provides the alternative representation

$$H(z,u) = s^{-t_0/s}(K(u)-z)^{-t_0/s}\mathcal{A}((1-u^h)(K(u)-z)^{h/s}), \tag{31}$$

which is our starting point. This expansion is analytically valid when (say) $|z| < \rho/2$ provided $u$ stays in a small enough neighborhood of 1. Write $\mathcal{A}(w) = \sum_{k\geq 0} a_k w^k$. One then has

$$H(z,u) = \sum_{k\geq 0} a_k (1-u^h)^k (K(u)-z)^{kh/s-t_0/s}. \tag{32}$$

Clearly, for the $r$th moment, it suffices to consider the sum in (32) with the index $k$ restricted to values in the interval $[0,r]$, so that

$$\chi_r(z) = \frac{1}{r!}\sum_{k=0}^{r} a_k (\partial_u^r((1-u^h)^k(K(u)-z)^{-t_0/s+kh/s}))_{u=1}. \tag{33}$$

The function $K(u)$ is analytic at $u=1$. Accordingly, the quantity $(K(u)-z)^{-1}$ and its derivatives at $u=1$ are of the form

$$(K(1)-z)^{-1} \equiv \frac{s}{1-sz}, \qquad -\frac{s^2K'(1)}{(1-sz)^2}, \qquad \frac{2s^3K'(1)^2}{(1-sz)^3} - \frac{s^2K''(1)}{(1-sz)^2},$$

and so on, with similar formulas holding for fractional powers. Thus $\chi_r(z)$ is invariably an algebraic function of a very special form, namely a finite linear combination of terms of the type

$$(1-sz)^{-t_0/s+kh/s-\ell}, \qquad 0 \leq k,\ell \leq r.$$

The statement then follows by coefficient extraction. □

As a consequence, one gets mechanically,

$$\mathbb{E}(X_n) \sim \frac{s+b}{s+h}sn, \qquad \mathbb{V}(X_n) \sim \frac{sh^2(s+a)(s+b)}{(s+h)^2(s+2h)}n,$$

which is consistent with the estimates of Bagchi and Pal [(1985), pages 395–397].

Finally, we turn to large deviations, for which the book of den Hollander (2000) can serve as a smooth introduction. It is known from the works of Hwang (1996) that a quasi-power approximation (in the sense of Lemma 1) for a family of PGFs leads to very precise "moderate deviation" estimates valid in some range not too far from the center of the distribution. We recycle here the technology of Hwang (1996), though the range is a little different. The large deviation rate is fully characterized by the following statement:

Corollary 4 (Large deviations). *Consider any balanced* $2 \times 2$ *urn with subtraction. Let* $\xi$ *be any number of the open interval* $(0, s\frac{s+b}{s+h})$*. One has*

$$\lim_{n\to\infty} \frac{1}{n} \log \mathbb{P}(X_n \le \xi \cdot n) = -\mathcal{R}(\xi), \tag{34}$$

*where the rate function* $\mathcal{R}$ *is determined from* $K(u)$ *defined in* (29) *by*

$$\mathcal{R}(\xi) = \max_{\lambda\in(0,1)} \log(s\lambda^{\xi} K(\lambda)). \tag{35}$$

*Equivalently,*

$$\mathcal{R}(\xi) = \log(s\lambda_0^{\xi} K(\lambda_0)) \quad \text{where } \lambda_0 \in (0,1) \text{ satisfies } \frac{\lambda_0 K'(\lambda_0)}{K(\lambda_0)} + \xi = 0. \tag{36}$$

[Put another way, $\mathcal{R}(\xi)$ is the Legendre transform of $\log(sK(e^t))$.]

Proof. Notice that $\mathbb{E}(X_n) \sim \frac{s+b}{s+h} sn$, so that (34) quantifies the left part of the distribution as approximately given by $e^{-n\mathcal{R}(\xi)}$. The basic ingredient is Cramér's technique of "shifting the mean" conjugated with upper bounds of the saddle point (equivalently, Chernoff) type as well as lower bounds based on the quasi-powers theorem in a shifted region.

First, one has

$$\mathbb{P}(X_n \le \xi n) = [u^k] \frac{p_n(u)}{1-u},$$

since multiplication by $(1-u)^{-1}$ sums coefficients of generating functions. Next, for any $f(u)$ analytic at 0 having nonnegative Taylor coefficients, the easy inequality $[u^k] f(u) \le f(\lambda)\lambda^{-k}$ holds provided the positive quantity $\lambda$ is taken inside the disc of convergence of $f(u)$. There results from these two observations the majorization

$$\mathbb{P}(X_n \le \xi n) \le \frac{p_n(\lambda)}{(1-\lambda)\lambda^{\lfloor \xi n \rfloor}}, \tag{37}$$

valid for any $\lambda \in (0,1)$.

In order to derive an *upper bound* on large deviations, it suffices to choose (as usual) the best possible value of $\lambda$ in (37). Now, for fixed positive $\lambda \in (0,1)$, the function $H(z,\lambda)$ has a dominant singularity of the algebraic type at $z = K(\lambda)$, see (23). A simple calculation based on the fact that the dominant singularities of $\psi$ are at $\rho\omega^j$ and that $I(u)$ increases from 0 to $I(1) = \rho$ for $u \in (0,1)$ shows further that $H(z,\lambda)$ has for $\lambda \in (0,1)$ a unique singularity at $K(\lambda)$ on $|z| = K(\lambda)$. Therefore one has by straight singularity analysis,

$$p_n(\lambda) \underset{n\to\infty}{\sim} C_\lambda \cdot s^{-n} K(\lambda)^{-n}, \tag{38}$$

for some constant $C$ (depending smoothly on $\lambda$). When $\xi$ lies in any fixed compact subinterval of $(0, s\frac{s+b}{s+h})$, the upper bound (37) can then be rewritten as

$$\mathbb{P}(X_n \le \xi n) \le \overline{C} s^{-n} K(\lambda)^{-n} \lambda^{-\xi n}$$

for some constant $\overline{C}$. This is a form amenable to optimization. Let $\lambda_0$ be such that $K(\lambda)^{-1}\lambda^{-\xi}$ attains its minimum over $(0,1)$ at $\lambda_0$. General convexity properties of probability generating functions imply that $\lambda_0$ exists and is unique.

The value of $\lambda_0$ is obtained by cancelling the derivative of $K(\lambda)^{-1}\lambda^{-\xi}$ and is thus a root of the second equation in (36). Up to factors that are subexponential in $n$, the upper bound in (38) is of the form $e^{-n\mathcal{R}(\xi)}$, with $\mathcal{R}(\xi)$ as given by (36) and (35). We have thus established "one half" of (34), namely,

$$\frac{1}{n} \log \mathbb{P}(X_n \le \xi n) \le -\mathcal{R}(\xi) + o(1),$$

with $\mathcal{R}(\xi)$ determined by (36).

There finally remains to argue that the upper bound is tight, that is, derive a *lower bound* on the probability values. This results from Cramér's technique of shifting the mean. The shifted law $r_{n,k} = [u^k] r_n(u)$ defined by the probability generating function

$$r_n(u) := \frac{p_n(\lambda_0 u)}{p_n(\lambda_0)}$$

satisfies a standard quasi-powers approximation and is itself amenable to Lemma 1. Assume first that the variability condition (26) holds for the shifted law given by $r_n(u)$. In that case the sum of probabilities $\sum_{\xi n - \sqrt{n} < k \le \xi n} r_{n,k}$ of the shifted law tends to a nonzero constant as it is approximated by a Gaussian integral. By construction, the $r_{n,k}$ are the $p_{n,k} \equiv [u^k] p_n(u)$ multiplied by a quantity $\lambda_0^k$ which varies between $e^{-O(\sqrt{n})}\lambda_0^{\xi n}$ and $O(1)\lambda_0^{\xi n}$. Thus, the corresponding sum $\sum_{\xi n - \sqrt{n} < k \le \xi n} p_{n,k}$ is, up to subexponential factors (themselves of the form $e^{-O(\sqrt{n})}$), of the type $e^{-n\mathcal{R}(\xi)}$. This implies a lower bound, hence the "other half" of the equality in (34). Finally, if the variability condition at $\lambda_0$ is not satisfied (this can only happen at isolated points), then an even stronger type of concentration holds for the shifted distribution $r_{n,k}$; in that case, the variance of the shifted distribution is $o(n)$, which, by Chebyshev's inequality, entails the stated lower bound on the sum of the $r_{n,k}$, hence the lower bound on partial sums of the $p_{n,k}$. □

The dual regime of large deviations on the right tail of the distribution is determined upon exchanging the roles of quantities $a$ and $b$.

**3. The $\mathcal{T}_{2,3}$ "tree" model.** The model is determined by its matrix and the initial conditions

$$\mathcal{T}_{2,3} = \begin{pmatrix} -2 & 3 \\ 4 & -3 \end{pmatrix}, \qquad a_0 = 2, \qquad b_0 = 0.$$

It has motivated much of the study of subtractive urns over the past two decades, given its relevance to several data structures of computer science [Aldous, Flannery and Palacios (1988), Bagchi and Pal (1985), Eisenbarth et al. (1982), Panholzer and Prodinger (1998) and Yao (1978)]. In this section, we arrive at the elliptic connection expressed by Theorem 3 and closely related to earlier works of Panholzer and Prodinger (1998). Our reason for treating this example in detail is twofold: first, it serves as a concrete illustration of the general treatment of Sections 1 and 2; second, it paves the way to our eventual characterization of the elliptic urn models in Section 4.

3.1. *Basic analytic structure.* Taking $a_0 = 2$ and $b_0 = 0$ corresponds to $t_0 = 2$. Theorem 1 provides an expression for $H(z,u)$:

$$H(z,u) = \delta(u)^2 \psi(z\delta(u) + I(u)).$$

Here $h = 6$, so that

$$\delta(u) := (1-u^6)^{1/6}, \qquad I(u) := \int_0^u \frac{t}{\delta(t)^5}\,dt = \int_0^u \frac{t}{(1-t^6)^{5/6}}\,dt,$$

and the function $\psi$ is defined implicitly by

$$\psi(I(u)) = J(u) \quad \text{where } J(u) := \frac{u^2}{\delta(u)^2} = \frac{u^2}{(1-u^6)^{1/3}}.$$

The results of Section 1 apply directly to this case. The elementary kite, which is the image of the sector $S_0$ of opening $\frac{\pi}{3}$, is a quadrilateral with vertices at $(0, \rho, I(\infty e^{2i\pi/12}), \rho\omega)$, where $\omega := e^{2i\pi/3}$. We find after a simple computation, upon following the proper branch of $\delta$,

$$\begin{aligned} I(+\infty e^{2i\pi/12}) &= \rho - e^{i\pi/6} \int_0^1 \frac{w^2\,dw}{(1-w^6)^{5/6}} \\ &= \rho - e^{i\pi/6}\frac{1}{6}B\left(\frac{1}{6}, \frac{1}{2}\right) = \rho\left(1 - e^{i\pi/6}\frac{\sqrt{3}}{2}\right). \end{aligned}$$

Thus, as $u$ varies from 0 to $+\infty$, passing through 1 (and above it), $I(u)$ describes first the segment from $[0,\rho]$, then the segment $[\rho, \rho(1+\omega)/2]$. The kite in this case happens to be a triangle and we shall refer to it as the "elementary triangle" (Figure 3).

There is also a "double parameterization" [due to evenness of both $I(u)$ and $J(u)$], so that we may freely identify points $u$ and $-u$. To this effect, we define

$$\mathcal{H} := \{z | (\Im(z) > 0) \vee ((\Im(z) = 0) \wedge (\Re(z) \geq 0))\}. \tag{39}$$

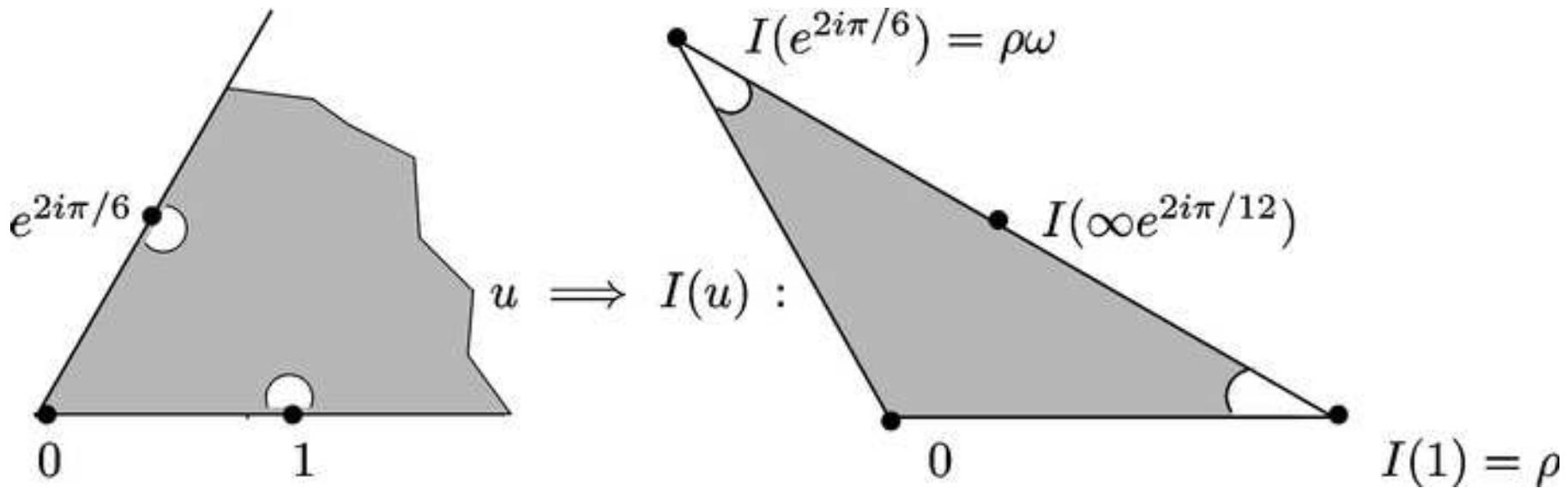

FIG. 3. *The "elementary kite," here a triangle,* $T_0$ (right) *is the image of the basic sector* $S_0$ (left) *via the mapping* $u \mapsto I(u)$.

Then, as $z$ ranges over $\mathcal{H}$, the fundamental polygon is obtained by gluing three rotated images of the elementary triangle. It is thus an equilateral triangle with center at the origin (Figure 4)—we call it the "fundamental triangle."

Next, the local analysis of $\psi$ at its dominant singularity $z = \rho$ results from the general treatment offered in Section 1. We find

$$\psi(z) = Z^{-2} - \tfrac{1}{7}Z^4 + \tfrac{1}{637}Z^{10} + \cdots, \qquad Z := \rho - z. \tag{40}$$

What is noteworthy here is the presence of a pole, rather than an algebraic singularity that prevails in the general case covered by Theorem 2. Similarly, the points $\rho\omega$ and $\rho\omega^2$ are double poles, so that the function

$$\psi(z) - \left(\frac{1}{(\rho - z)^2} + \frac{1}{(\rho\omega - z)^2} + \frac{1}{(\rho\omega^2 - z)^2}\right) \tag{41}$$

is analytic in a disc $|z| < R$ for some $R > \rho$.

The fact that the dominant singularities of $\psi$ are poles naturally led us to look for the next layer of singularities, as this would provide very precise information on the exponential smallness of error terms. In so doing, much to our surprise, we uncovered a lattice structure commonly associated with elliptic functions (Figure 5).

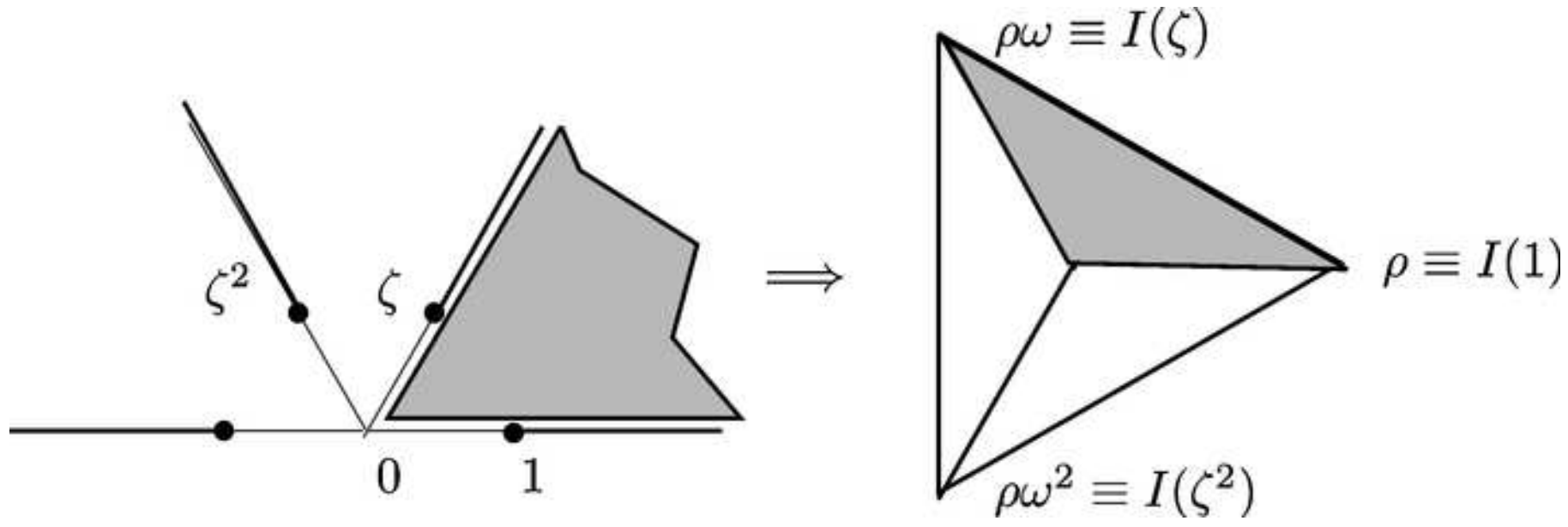

FIG. 4. *The "fundamental polygon," here a triangle,* $T$ (right) *is the image of the slit upper half-plane* ($\mathcal{H}$) (left) *via the mapping* $u \mapsto I(u)$.

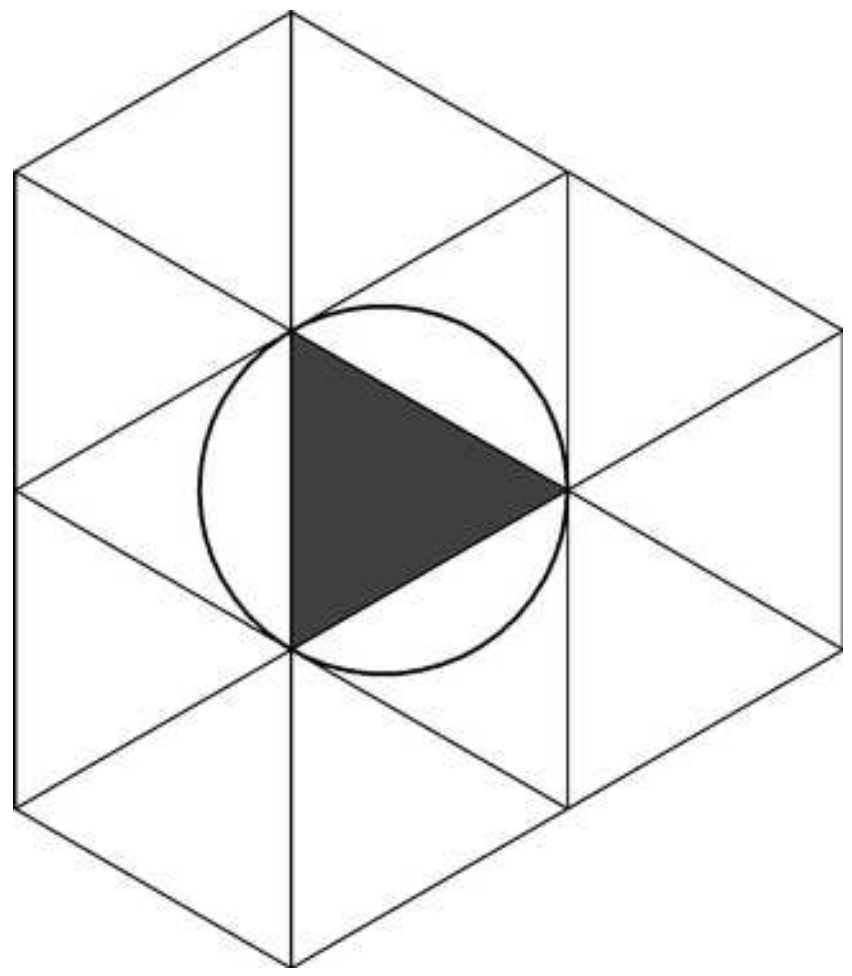

FIG. 5. *Rotated copies of the fundamental triangle around $\rho, \rho\omega, \rho\omega^2$ shown against the circle of convergence of $\psi(z)$.*

3.2. *The elliptic structure.* An *elliptic* function is a function that is meromorphic in the whole complex plane and is doubly periodic. Amongst the many different ways to develop the corresponding theory, perhaps the simplest is the one originally proposed by Weierstraß, where elliptic functions are defined as sums of rational functions taken over lattices. [Accessible introductions appear in the books by Whittaker and Watson (1927) and Chandrasekharan (1985).]

DEFINITION 2. A *lattice* $\Lambda$ with generators $\xi, \eta \in \mathbb{C}$ is defined as the set of complex numbers

$$\Lambda(\xi, \eta) = \{n_1\xi + n_2\eta | n_1, n_2 \in \mathbb{Z}\}.$$

The Weierstraß $\wp$-*function* relative to $\Lambda$ is classically defined as

$$\wp(z; \Lambda) = \frac{1}{z^2} + \sum_{w \in \Lambda \setminus \{0\}} \left( \frac{1}{(z-w)^2} - \frac{1}{w^2} \right). \tag{42}$$

(The Weierstrass $\wp$-function is by construction doubly periodic.)

We shall make use here of the "hexagonal" lattice $\Lambda$ defined as the lattice generated by $e^{i\pi/6}, e^{-i\pi/6}$, see Figure 5,

$$\Lambda_{\text{hex}} := \{n_1 e^{i\pi/6} + n_2 e^{-i\pi/6} | n_1, n_2 \in \mathbb{Z}\}, \tag{43}$$

and its associated Weierstraß zeta function, $\wp(z; \Lambda_{\text{hex}})$. We state:

Theorem 3 (Elliptic connection). *The $\psi$-function of the $\mathcal{T}_{2,3}$ model initialized with two balls of the first type ($a_0 = t_0 = 2$) is exactly*

$$\psi(z) = \frac{1}{(\rho\sqrt{3})^2}\wp\left(\frac{z-\rho}{\rho\sqrt{3}}\right) \qquad \text{with } \rho := \frac{1}{6}\frac{\Gamma(1/3)\Gamma(1/6)}{\Gamma(1/2)}, \tag{44}$$

*where $\wp(z) := \wp(z; \Lambda_{\text{hex}})$ is the Weierstraß function of the hexagonal lattice. In particular, the bivariate generating function of the model is expressible in terms of elliptic functions.*

Note. The function $\psi$ can be alternatively written as $\psi(z) = \wp(z - \rho \,|\, 0, -4)$ where $\wp$ is specified by the lattice invariants $g_2 = 0$ and $g_3 = -4$.

Proof of Theorem 3. Consider the whole complex plane tiled by nonoverlapping copies of the hexagon of center $\rho$, radius $\rho\sqrt{3}$, having vertices at the points $\rho + \rho\sqrt{3}\Lambda_{\text{hex}}$.

We claim that any complex point $z$ is reachable as a value $I(\gamma(u))$, where the notation $I(\gamma(u))$ indicates that the integral defining $I$ is to be taken along a path $\gamma(u)$ that starts at 0 and ends at $u$. Similarly, $J(\gamma(u))$ will represent the determination of $J(u)$ along path $\gamma(u)$ that is obtained by continuity from the principal determination at 0. Otherwise said, we are walking on the Riemann surface of the Fermat curve $\delta(u)$.

The algorithm is as follows. Assume for simplicity that $z$ is the center of one of the equilateral triangles in which the hexagonal tiling decomposes. The straight line $L_0$ from 0 to $z$ can be first slightly deformed into a curve $L_1$ that avoids all the vertices of the tiling. This $L_1$ can then be transformed into a polygonal line $L_2$ that connects centers of successive equilateral triangles. Finally, each segment of $L_2$ can be changed into a pair of segments going through one of the vertices of the lattice and forming an angle a positive multiple of $\pi/3$. The resulting polygonal line, $L_3$, will be called the *standard $z$-path.* See Figure 6 for a graphic rendering.

The contour $\gamma$, called the *standard $u$-path*, is then obtained from the standard $z$-path $L_3$ by first applying a contraction by a factor $\rho$, then executing the following routine:

- turn by an angle of $6\theta$ whenever $L_3$ turns at an angle of $\theta$ (where $\theta$ is a multiple of $\pi/3$) around a vertex of the lattice,
- turn by an angle $\theta/2$ (where $\theta$ is a multiple of $2\pi/3$) whenever $L_3$ turns by $\theta$ around the center of one of the equilateral triangles.

The construction is then easily modified to accommodate points that are not centers of triangles of the tiling.

For any $z$ in $\mathbb{C}$ that is not a vertex of the tiling, the algorithm described above determines constructively a path $\gamma(u)$. By design, along such a path, one has $I(\gamma(u)) = z$. Indeed, the standard $u$-path is precisely such that it "undoes" the effect of $I(u)$ on angles at points either vertices of the tiling or centers of the triangles; at the same time, the variation of $I(u)$ along a segment from a point $u_0$ above 0 to a point $u_1$ above some $\zeta^j$ is precisely of modulus $\rho$ and thus gives rise to a segment with the "right" length. See once more Figure 6. In this way, we find that $I(\gamma(u))$ reaches any point $z$ of the complex plane that is not a vertex of the tiling, and at the final point, $J(\gamma(u))$ is locally analytic, so that $\psi$ is itself analytic at $z$. Thus, $\psi(z)$ can be continued to the complex plane punctured at vertices of the tiling.

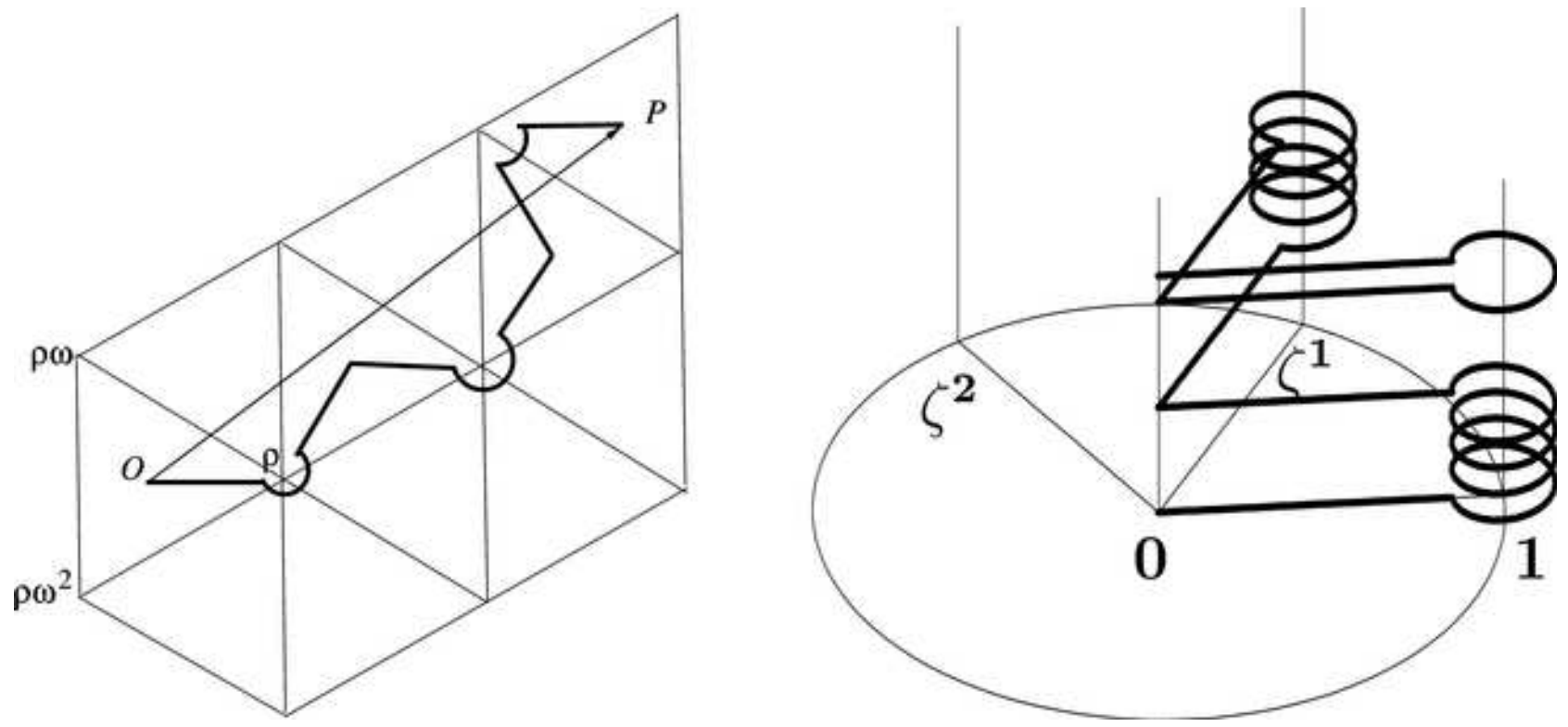

FIG. 6. *A standard path in the $z$-plane from* 0 *to* $P \equiv z$ *and the contour* $\gamma$ *above the* $u$*-plane that realizes it via* $u \mapsto z = I(\gamma(u))$.

When $z = w$ is one vertex of the lattice, then it is approached from a certain direction by a path $\gamma(u)$, where $u$ is near $\zeta^0$, $\zeta^1$ or $\zeta^2$. Along the path a certain determination $\delta^\circ(u)$ of $\delta$ is in force, where all determinations are of the form $\zeta^r \delta(u)$ with $0 \le r < 6$. Then, the very same determination $\delta^\circ$ must be adopted in $J(\gamma(u))$ that tends to infinity as $I(\gamma(u))$ approaches $w$. A local analysis entirely analogous to the one conducted for the three dominant poles shows that $\psi$ has a double pole at $w$, and that its principal part there consistently exhibits the same dominant coefficient and residue.

The analytic continuation of $\psi(z)$ along such paths $\gamma$ therefore has the same dominant parts and residues at double poles as the right-hand side of (44), namely the function

$$\widetilde{\psi}(z) := \frac{1}{(\rho\sqrt{3}\,)^2}\wp\left(\frac{z-\rho}{\rho\sqrt{3}}\right).$$

Consequently, the difference $\psi(z) - \widetilde{\psi}(z)$ is an entire function. That this entire function reduces to 0 results from Liouville's theorem, as we finally argue.

Draw discs of some sufficiently small but fixed radius around the six roots of unity in the $u$-plane and consider these as excluded regions in the con-struction of $u$-paths. Then the image $z = I(\gamma(u))$, as $\gamma(u)$ varies, avoids the plan stripped of small ovals around the corresponding lattice points of the $z$-plane. But, at the same time $J(\gamma(u))$ remains bounded. Therefore, on the complex plane with "holes," $\psi(z)$ is uniformly bounded by a constant. From the fact that $\widetilde{\psi}$ is doubly periodic, there results that it is also bounded over the plane with holes, hence

$$|\psi(z) - \widetilde{\psi}(z)| < c_1,$$

for some $c_1 > 0$. In particular, the bound holds on an infinity of near-circular contours centered at the origin and having arbitrarily large diameter. Then, by virtue of a known variant of Liouville's theorem (an entire function bounded in modulus along large contours is a constant), one must have identically

$$\psi(z) - \widetilde{\psi}(z) = d_1,$$

for some complex constant $d_1$. This constant is actually equal to 0 as is seen from comparing the expansions of $\psi(z)$ and $\widetilde{\psi}(z)$ at 0. The proof of the theorem is completed. □

3.3. *Probabilistic consequences of the analytic model for $\mathcal{T}_{2,3}$.* We are now in a position to exploit the analytic solutions expressed by Theorem 3. The general theory of Sections 1 and 2 applies, giving the large deviation rate function and the limit Gaussian law. In addition, curious exact representations as sums over lattice points result for the probability generating functions describing the urn composition (Section 3.3.1). Surprisingly perhaps, a very precise form of all moments can be obtained in terms of a family of polynomials of "binomial type" [Rota (1975)]; see Section 3.3.2.

3.3.1. *Exact representations and Gaussian laws.* The lattice structure that underlies the Weierstraß function is directly reflected at the level of coefficients. The resulting form below is naturally very strong, as it is an *exact* description of the probability generating function at time $n$.

COROLLARY 5 (Elliptic structure of $\mathcal{T}_{2,3}$). *For the $\mathcal{T}_{2,3}$ model, the probability generating function $p_n(u) = \mathbb{E}(u^{X_n})$ admits an* exact *formula valid for all $n \geq 1$,*

$$p_n(u) = \sum_{n_1,n_2=-\infty}^{+\infty} \left( K(u) + \frac{\rho\sqrt{3}}{\delta(u)} (n_1 e^{i\pi/6} + n_2 e^{-i\pi/6}) \right)^{-n-2}, \tag{45}$$

*where*

$$K(u) := \frac{1}{\delta(u)} \int_u^1 \frac{t}{\delta(t)^5}\, dt, \qquad \delta(u) = (1-u^6)^{1/6}.$$

PROOF. From Theorem 3, we need to extract $[z^n]\delta^2\psi(\delta z + I)$, where $\psi$ admits a decomposition as a sum of rational fractions over elements of the lattice $\Lambda$. Then, after a simple calculation, one gets

$$p_n(u) = \frac{1}{n+1} \sum_{w \in \Lambda^\star} \left( [z^n] \frac{1}{(w-z)^2} \right),$$

where $\Lambda^\star$ is a translated and scaled version of $\Lambda$:

$$\Lambda^\star = \frac{\rho\sqrt{3}}{\delta(u)}\Lambda + K(u).$$

The result follows. □

Corollary 2 regarding general Gaussian limits applies here. Thus, for the $\mathcal{T}_{2,3}$ model, the random variable $X_n$ representing the number of balls of the first type at time $n$ is asymptotically Gaussian with speed of convergence to the limit $O(n^{-1/2})$, in the sense of (27). The random variable $X_n$ superficially resembles a sum of independent random variables since its probability generating function is essentially an $n$th power of the fixed function $K(u)^{-1}$. It is, however, of interest to observe that the function $K(u)^{-1}$, though analytic at 0 and satisfying $K(1) = 1$, is *not* a probability generating function, as its Taylor coefficients of index $6, 12, 18, \ldots$ turn out to be negative:

$$K(u)^{-1} \doteq 0.713 + 0.254u^2 + 0.090u^4 - 0.086u^6 + 0.022u^8 + \cdots.$$

3.3.2. *The shape of moments.* An interesting consequence of the elliptic connection concerns moments of the distribution of the urn's composition, $X_n$. Bagchi and Pal, Mahmoud, and Panholzer and Prodinger have determined the exact form of the first two moments, while Bagchi and Pal have obtained further asymptotic information on the moments of higher order. This already involved a certain amount of calculational effort with recurrences. In fact, globally, the moments have an amazingly simple form deriving from the elliptic connection.

COROLLARY 6 (Moments of $\mathcal{T}_{2,3}$). *For the $\mathcal{T}_{2,3}$ model, exact polynomial forms of moments of* any *order are available: the factorial moments satisfy*

$$\mathbb{E}((X_n)^{\underline{r}}) = P_r(n+2), \qquad n \geq 6r - 1,$$

*where the $P_r$ are polynomials generated by*

$$e^{vL(h)} = \sum_{r=0}^{\infty} \frac{h^r}{r!} P_r(v) \quad \text{and} \quad L(h) = -\log K(1+h). \tag{46}$$

Using a symbolic manipulation system, the polynomials are easily computed from the expansion of $K$ at $u=1$. To wit:

$$K(1+h)=1-\tfrac{4}{7}h+\tfrac{10}{91}h^2+\tfrac{300}{1729}h^3-\tfrac{1689}{8645}h^4+\cdots.$$

One then finds mechanically

$$P_1(\nu)=\frac{4\nu}{7},\qquad P_2(\nu)=\frac{4\nu}{637}(52\nu+17),$$
$$P_3(\nu)=\frac{8\nu}{84721}(1976\nu^2+1938\nu-11063).$$

In particular, the mean and variance of $X_n$ are

$$\mathbb{E}(X_n)=\tfrac{4}{7}(n+2),\qquad \mathbb{V}(X_n)=\tfrac{432}{637}(n+2)^2.$$

PROOF OF COROLLARY 6. Take the fundamental PDE, isolate $G'_u(z,u)$ and repeatedly differentiate with respect to $u$, then set $u=1$. This provides a triangular system from which one can "pump" in succession the generating functions of moments of order $1,2,3,\ldots.$ One then verifies by induction that the ordinary generating function of the moments of order $r$ is of the form

$$\sum_n \mathbb{E}(X_n^{\underline{r}})z^n=\frac{\widetilde{P}_r(z)}{(1-z)^{r+1}}+\widetilde{Q}_r(z),$$

where $\widetilde{P}_r,\widetilde{Q}_r$ are polynomials and

$$\deg(\widetilde{P}_r(z))\le r,\qquad \deg(\widetilde{Q}_r(z))\le 6r-2.$$

This argument grants us *nonconstructively* the existence of a polynomial representation for each moment as soon as $n$ is large enough.

There remains to identify the particular class of polynomials involved. Start from the fact that

$$p_n(u)=K(u)^{-n-2}+\text{exponentially small terms in } n.$$

Since the factorial moment of order $r$ satisfies

$$\mathbb{E}(X_n^{\underline{r}})=(\partial_u^r p_n(u))_{u=1}=[(u-1)^r]p_n(u),$$

it can be obtained, up to exponentially small error terms, by expanding $K(u)^{-n-2}$ around $u=1$. Retaining only the polynomial part (in $n$),

$$[(u-1)^r]K(u)^{-n-2}=[(u-1)^r]e^{-(n+2)\log K(u)}=[h^r]e^{-(n+2)\log K(1+h)},$$

we get what the statement asserts. □

3.3.3. *Large deviations.* Corollary 4 applies to the effect that the large deviation rate function is a transform of $K(u)$. An immediate consequence of the analysis of the polar singularities of $\psi$ is a very precise quantification of *extreme large deviations*:

COROLLARY 7 (Extreme large deviations of $\mathcal{T}_{2,3}$). *The probability that, at time* $3n+1$ *in the* $\mathcal{T}_{2,3}$ *model, all balls are of the second color* (W) *is (any* $A<8$*):*

$$3\rho^{-3n-3}(1+O(A^{-n})).$$

PROOF. The function $\psi(z)$ is exactly the BGF of the urn at $u=0$. Thus $(n+1)^{-1}[z^n]\psi(z)$ is the probability for the urn not to contain any ball of the first type. The property then results immediately from the fact that the next layer of poles is on $|z|=2\rho$. □

**4. The elliptic cases.** In this section, we list all the cases of urns with subtraction which, like the $\mathcal{T}_{2,3}$ model, lead to elliptic functions; such urns are, as in the previous section, attached to exact lattice representations and simplified moment forms. We shall say that an urn is *elliptic* if, for some choice of initial conditions $(a_0,b_0)$, the fundamental function $\psi$ is a fractional power of an elliptic function,

$$\psi(z)=\Pi(z;\Lambda)^{p/q}, \tag{47}$$

where $\Pi$ is meromorphic and admits a period lattice $\Lambda$. The power will depend on the initial conditions, but by taking the initial configuration determined by $a_0$ and $b_0$ sufficiently large, one can always render the exponent integral; see the remarks on "sensitivity to initial conditions" as well as (19) and (20). Obviously, for urns with subtraction, one need only consider models that are *arithmetically irreducible*: by this is meant that the matrix $M=\binom{\alpha\ \ \beta}{\gamma\ \ \delta}$ defining the urn has coprime entries: $\gcd(\alpha,\beta,\gamma,\delta)=1$.

The key characters in this section are the following six urns:

$$\begin{aligned} &A=\begin{pmatrix}-2 & 3\\ 4 & -3\end{pmatrix}, \qquad B=\begin{pmatrix}-1 & 2\\ 3 & -2\end{pmatrix}, \qquad C=\begin{pmatrix}-1 & 2\\ 2 & -1\end{pmatrix},\\ &D=\begin{pmatrix}-1 & 3\\ 3 & -1\end{pmatrix}, \qquad E=\begin{pmatrix}-1 & 3\\ 5 & -3\end{pmatrix}, \qquad F=\begin{pmatrix}-1 & 4\\ 5 & -2\end{pmatrix}.\end{aligned} \tag{48}$$

The urn $A$ is of course the $\mathcal{T}_{2,3}$ model. As is easily checked, any urn of balance $s=1$ is necessarily of type $A$, $B$ or $C$ and any arithmetically irreducible urn of balance 2 can only be of type $D$ or $E$. Thus, the first five cases exhaust all possible types of urns with subtraction having balance $s=1,2$. The urn $F$ is one of the four possible irreducible urns having balance $s=3$. We state:

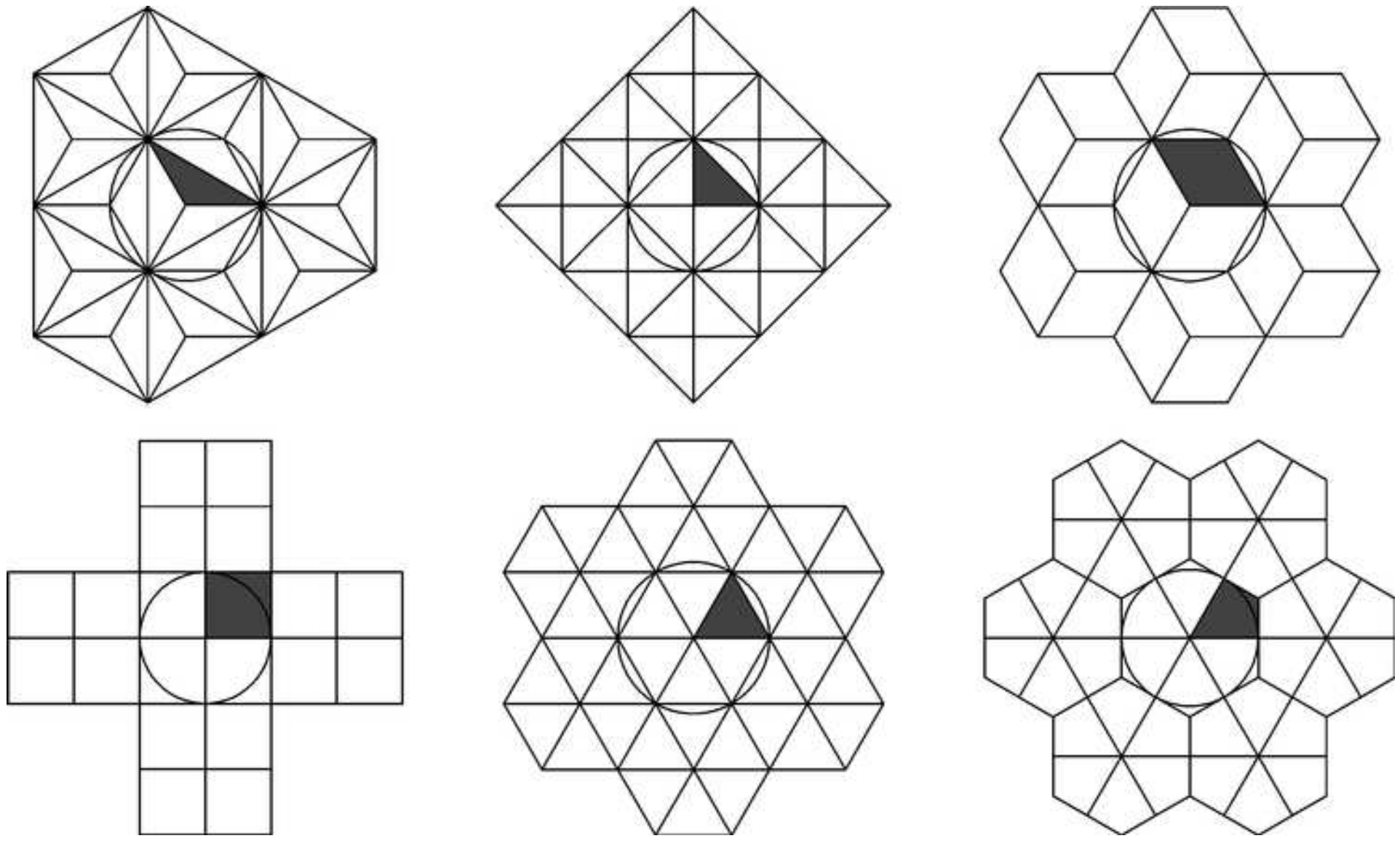

FIG. 7. *The six elliptic cases in order $A, B, C, D, E, F$: the diagrams formed by the fundamental polygon together with its rotated images. (The elementary kite is darkened.)*

THEOREM 4. *All balanced $2 \times 2$ urns with balance $s = 1$ (cases $A, B$ and $C$) are elliptic. All urns with balance $s = 2$ (cases $D, E$) are also elliptic. There exists one "sporadic" urn (case $F$) of balance $s = 3$ that is elliptic. These six matrices represent the only arithmetically irreducible models of urns with subtraction that are elliptic.*

PROOF. For 2–3 trees, the matrix $A$ corresponds to $a = 2, b = 3, s = 1, h = 6$. As seen in Section 3.2, it is associated to a regular tiling of the plane by equilateral triangles. Its $\psi$ function has, by virtue of Theorems 2 and 3, a *double pole* at $\rho$ and all other points of the lattice and is a Weierstraß $\wp$-function.

Next, we turn to the other urns of balance 1. The corresponding mapping properties are represented in Figure 7. The model $B$ has an elementary kite which is a right triangle with vertices $0, \rho, i\rho$, so that the fundamental polygon is the square with vertices $\rho, i\rho, -\rho, -i\rho$. This fundamental square tiles the plane. The model $C$ leads to a lattice much similar to the 2–3 tree case. Theorem 2 shows that in all the three cases $A, B$ and $C$, the $\psi$ function has a pole at $\rho$ since the principal exponent at the singularity $t_0/s$ as well as the Puiseux exponent $h/s$ of (23) are integers (here $s = 1$). Double periodicity then results from arguments similar to those developed in the proof of Theorem 3.

An analogous discussion applies to the two urns of balance 2, namely $D$ and $E$. In the case of $D$, the function $I(u)$ is even exactly a lemniscatic

integral:

$$I(u) = \int_0^u \frac{dt}{\sqrt{1-t^4}}.$$

Finally, among the four urns of balance 3, namely

$$\begin{pmatrix} -1 & 4 \\ 4 & -1 \end{pmatrix}, \qquad \begin{pmatrix} -1 & 4 \\ 5 & -2 \end{pmatrix}, \qquad \begin{pmatrix} -1 & 4 \\ 7 & -4 \end{pmatrix}, \qquad \begin{pmatrix} -2 & 5 \\ 8 & -5 \end{pmatrix},$$

one found to be elliptic has $s = 3$, $h = 6$, so that $h/s = 2$ is integral: this is $F = \left(\begin{smallmatrix} -1 & 4 \\ 5 & -2 \end{smallmatrix}\right)$. By the usual reasoning, the function $\psi$ leads to an elliptic function when one starts with $a_0 = 3$ and $b_0 = 0$.

Note a necessary condition for an urn to be elliptic: in addition to the tenability conditions, $s$ should divide $h$, in accordance with Theorem 2 and (23); that is, the Puiseux exponent $h/s$ should be integral. (Otherwise, all powers of $\psi$ inherently have branch points and hence cannot be meromorphic.) The other three cases of balance $s = 3$, including the "pentagonal" urn of Figure 2, correspond to a fractional value of the Puiseux exponent $h/s$ in (23) and therefore cannot be reduced to elliptic functions.

There finally remains to prove that all elliptic urns have indeed been found. The condition that the Puiseux exponent $h/s$ is integral, arithmetic irreducibility, and the tenability conditions taken together imply the existence of triples $(x, y, z)$ of integers, representing up to possible permutation the values $(a, b, s)$, such that

$$\gcd(x, y, z) = 1, \quad x \mid y + z, \quad y \mid z + x, \quad z \mid x + y. \tag{49}$$

Simple arithmetic shows that the only values for which the system admits a solution are permutations of the basic types

$$(1, 1, 1), \quad (1, 1, 2), \quad (1, 2, 3), \tag{50}$$

and, in particular, one must have $s \leq 3$. The arithmetic argument goes as follows. Take $x \leq y \leq z$. One has $z \leq 2y$ since $z \leq x + y$ by the third divisibility condition in (49); then note the stronger property that $x, y, z$ are pairwise coprime (proof: a contrario). Then set $z + x = qy$, where $q \leq 3$; the first divisibility condition then implies $x|(q + 1)y$, and since $x$ and $y$ are coprime, $x|4$ so that $x \in \{1, 2, 4\}$. This in turn implies $y \leq z \leq y + 4$ while $z$ must divide $y + x$. Combining this with the second and third divisibility conditions of (49), we see that there are only finitely many possibilities which are then easily tested. Finally, completeness of the list (50) implies that elliptic cases have indeed all been found. □

It is pleasant to note that the elliptic urns correspond to the *crystallographic groups* of the Euclidean plane, that is, groups of isometries acting discontinuously (in fact, the ones which admit a compact fundamental domain). As is well known, these groups themselves describe the possible regular tessellations of the plane by polygonal tiles and are in finite number; see, for example, Berger (1977) and Yoshida (1997).

**5. Discussion.** A probabilist may legitimately expect more than standard-issue central limit theorems and Cramér approximations to come out of an analytic treatment of urn models. We would like to offer a few comments.

1. The combinatorial derivation to the fundamental PDE of (14) applies (with only notational adjustments) to any urn model with *more than two colors*, provided it remains *balanced*. The resulting PDE is invariably first order and linear. This study has exhibited a class, the subtractive balanced $2 \times 2$ models, for which the associated ordinary differential system provided by the method of characteristics and the corresponding generating functions prove to be analytically tractable. This suggests to look for other cases where analysis can be made to work.

2. On the algebraic–analytic front, Theorems 1 and 2, though they appear to admit no universal $r \times r$ extension when $r > 2$, can at least be adapted and generalized to several models [work in progress with V. Puyhaubert (2003–2004)]:

(i) $2 \times 2$ balanced urns with *positive* entries of the type of "Bernard Friedman's urn" [Friedman (1949)]. Interestingly enough, an adapted form of the algebraic solution expressed by Theorem 1 holds (only integration constants need to be changed). This may be seen as an elicitation of some of Friedman's remarks concerning his differential recurrences [Friedman (1949), pages 61 and 62]. Developments parallel to Theorem 2 appear to be possible, they then open access to *non-Gaussian laws* which have not been previously made explicit.

(ii) $2 \times 2$ balanced urns corresponding to *nonnegative triangular* matrices. Gouet (1993) shows strong functional limits to exist, but some of the involved characteristics have remained inaccessible due to nonconstructive aspects of martingale theory. In that case, we can supplement Gouet's finding and derive explicitly *stable laws* and a *local limit* theorem.

(iii) Balanced urns with *three colors* and a nonnegative *triangular* matrix, as well as some special cases of triangular $r \times r$ urns for $r > 3$.

3. The determination of the *large deviation rate* is obtained by standard arguments of Cramér type, once the relevant analytic forms have been established, but, to the best of our knowledge, it is new. We observe that the asymptotic approximations we have obtained, when combined with the saddle point method, can provide strong forms of large deviation estimates, complete with subexponential factors and multipliers. As we have seen in the case of moment estimates and Gaussian laws, a precise determination of the speed of convergence to the asymptotic regime is a boon provided by the analytic machinery. We also observe that the general algebraic solutions supplied by Theorem 1, when coupled with (19) and (20) (sensitivity to initial conditions), make it possible, in principle, to analyze the evolution of the urn starting with a large number of balls of each color. This could be put to use in order to analyze finely the distribution of sample paths, in the central or large deviation regime.

4. Besides integer partitions and theta functions, *elliptic models* occasionally pop up in combinatorics, some models related to permutations having been found by Dumont, Flajolet, Françon and Viennot about 1980. The present work contributes a new kind, *urn histories*, themselves equivalent to certain weighted lattice paths. On another register, given that algebraic functions of higher genus can be uniformized by theta-Fuchsian functions, we may even fantasize about the possibility of "nonprobabilistic" representations for urns of genus higher than 1 in terms of tessellations of the hyperbolic plane.

**Acknowledgments.** We are grateful to Philippe Robert, Hosam Mahmoud, Frédéric Chazal, Michel Merle and Arnaud Beauville for stimulating discussions and bibliographic data. We would also like to thank an anonymous referee for suggesting a major reorganization of an earlier version, resulting, we hope, in a clearer presentation of our work.

P. FLAJOLET
PROJECT ALGORITHMS
INRIA ROCQUENCOURT
78153 LE CHESNAY
FRANCE
E-MAIL: Philippe.Flajolet@inria.fr

J. GABARRÓ
H. PEKARI
DEPARTAMENT DE LLENGUATGES
I SISTEMES INFORMÀTICS
UNIVERSITAT POLITÈCNICA DE CATALUNYA
CAMPUS NORD, EDIF. C5
JORDI GIRONA 1-3
E-08034 BARCELONA
SPAIN
E-MAIL: gabarro@lsi.upc.es
E-MAIL: lsi@pekari.de